\documentclass[fleqn,sort&compress,twocolumn]{svjour3mod}

\usepackage{graphicx}
\usepackage{bm}
\usepackage{amsmath,amsfonts,cuted} 
\usepackage{booktabs}
\usepackage{multirow}

\newcounter{myremark}
\newenvironment{myremark}[1][]{\refstepcounter{myremark}\par\medskip
   \noindent \textbf{Remark~\themyremark.
	 #1} \rmfamily}{\medskip}

\usepackage{natbib}
\setcitestyle{numbers,square}
\usepackage[svgnames]{xcolor}
\usepackage[colorlinks]{hyperref}
\AtBeginDocument{%
  \hypersetup{
    citecolor=Indigo,
    linkcolor=DarkGreen,   
    urlcolor=Gray}}

\newif\ifdraftversion
\draftversionfalse

\ifdraftversion
  \usepackage[colorinlistoftodos, obeyFinal]{todonotes}
  \newcommand\pablo[1]{\todo[inline,color=red!50!white]{Pablo commented: #1}}
  \newcommand\thibaut[1]{\todo[inline,color=blue!50!white]{Thibaut commented: #1}}
\else
  \newcommand\pablo[1]{}
  \newcommand\thibaut[1]{}
\fi

\def\beme{{\bm{m}}}
\def\bn{{\bm{n}}}

\def\bN{{\bm{N}}}
\def\bM{{\bm{M}}}

\def\bRhi{{\hat{\bm{R}}_i}}
\def\rhi{{\hat{r}_i}}

\def\diff{{\text{d}}}
\def\bx{{\bm{x}}}
\def\bxh{{\hat{\bx}}}

\def\xt{{\tilde{x}}}
\def\dbx{{\diff\bm{x}}}
\def\dbxh{{\diff\hat{\bx}}}

\def\dxt{{\diff\tilde{x}}}

\def\Si{{\bm{S}_i}}

\def\Fi{{F_i}}
\def\Fhi{{\hat{F}_i}}

\def\ghij{{\hat{\gamma}_{i,j}}}

\def\chij{{\hat{\bm{c}}_{i,j}}}

\def\hnabla{{\hat{\nabla}}}

\def\T{{\mathcal{T}}}

\def\sumi{{\sum_{i=1}^{n_F}}}
\def\sumj{{\sum_{j=1}^{n_{c,i}}}}

\newcommand{\bA}{\ensuremath{\bm{A}}}

\newcommand{\diffusivityCoef}{\bm{\ensuremath{K}}}

\newcommand{\Vi}[2]{\bm{#1}_{#2}}

\begin{document}\sloppy

\journalname{arXiv math.NA}
\title{Quadrature-free Immersed Isogeometric Analysis}

\author{P. Antolin \and T. Hirschler}
\institute{P. Antolin - T. Hirschler \at
  Institute of Mathematics, Chair of Numerical Modelling and Simulation\\
  École Polytechnique Fédérale de Lausanne, Switzerland 
  \email{pablo.antolin@epfl.ch}
}
\date{}

\maketitle

\begin{abstract}
    This paper presents a novel method for solving partial differential equations on three-dimensional CAD geometries by means of immersed isogeometric discretizations that do not require quadrature schemes.
    It relies on a new developed technique for the evaluation of polynomial integrals over spline boundary representations that is exclusively based on analytical computations.
    First, through a consistent polynomial approximation step, the finite element operators of the Galerkin method are transformed into integrals involving only polynomial integrands.
    Then, by successive applications of the divergence theorem, those integrals over B-Reps are transformed into first surface and then line integrals with polynomials integrands. Eventually these line integrals are evaluated analytically with machine precision accuracy.
    The performance of the proposed method is demonstrated by means of numerical experiments in the context of 2D and 3D elliptic problems, retrieving optimal error convergence order in all cases.
    Finally, the methodology is illustrated for 3D CAD models with an industrial level of complexity.

\end{abstract}

\section{Introduction}
\label{sec:intro}

The integration of Computer-Aided Design (CAD) and Computer-Aided Engineering has gained interest during the last two decades with the introduction of new numerical approaches as, for instance, the isogeometric paradigm~\cite{Hughes_2005,Cottrell_2009} or meshfree strategies~\cite{Liu_2009}.
Particularly, spline-based geometric models have been found to present excellent performance for numerical simulations~\cite{BAZILEVS_2006,Buffa_2011,Hiemstra_2014,Lipton_2010}.
This opens the door to the formation of all-in-one design frameworks where a single geometric model is simultaneously used for parameterizing the shape of the object of interest and performing advanced numerical analyses~\cite{Herrema_2017,Antolin_2019,Hafner_2019,Hirschler_2020}.
The combination into one single model of both high-fidelity geometrical properties and efficient analysis performances is however far from trivial in general.
Indeed, generating analysis-suitable geometric models for complex industrial designs requires advance numerical tools.
To achieve this goal, two different strategies can be undertaken: The first one consists in generating a fully conformal multi-patch geometric model such that standard analysis procedures can be directly employed.
Generating these conformal meshes is however a quite challenging task in the case of geometries with complex topologies~\cite{wang2012converting,wei2018blended,xia2017isogeometric}, especially when only tensor-product splines are considered~\cite{Al_Akhras_2016,Hinz_2018,Massarwi_2019,Maquart_2020}.
On the contrary, the second approach aims to directly use standard CAD models which may contain non-conforming and trimmed surfaces and present geometric defects, as water leaks or surface overlaps, and to recall to high-end analysis procedures~\cite{rank2012geometric,Legrain_2013,Breitenberger_2015,Hsu_2016,guo2018variationally,wassermann2019integrating}.
Interest readers may refer to~\cite{Marussig_2017}, and the many references therein, for an extensive review in the context of isogeometric methods.
The present work falls into this second category.

A major ingredient that is commonly required in order to perform numerical analyses over CAD models is an efficient integration procedure which enables to evaluate integrals over complex domains such as curved polyhedrons.
This is, for instance, the case when employing non-conformal analysis methods, where the geometric representation is decoupled from the discretization of the solution~\cite{peskin2002immersed,D_ster_2008,schillinger2012isogeometric,burman2015cutfem,wassermann2017geometric,elfverson2018cutiga}.

In this context of immersed and enriched FEM, there exist several integration approaches.
In 3D, among the most common ones is worth highlighting \emph{octree subdivision}~\cite{shephard1991automatic,abedian2013performance,kudela2016smart,peto2020enhanced} which consists in adaptively subdividing the domain of integration into sub-cells (voxels in 3D, or simple pixels in 2D).
The obtained piecewise constant approximation of the underlying geometry can be improved by performing a local boundary reparameterization at the finest level of this recursion procedure via a (low-order) tessellation method~\cite{Verhoosel_2015,Divi_2020}.
Despite the beneficial simplicity and robustness of this decomposition-based method, it may suffers from high computational cost due to the large number of integration sub-cells, especially in three-dimensional and high-order methods.

For problems where the geometric representation of the boundary is of major importance, alternative approaches are considered as for instance \emph{high-order triangulation} and other sophisticated partitions~\cite{Kudela_2015,kudela2016smart,Massarwi_2019,Antolin_2019b}.
They consist in generating boundary-conforming sub-meshes which are generally non analysis-suitable (due to the presence of hanging nodes, missing connectivity, singularities, etc.) but which are handy for integration purposes.
The high-fidelity representation of the geometry boundaries, even for complex geometries, yields a high-accuracy in the evaluation of integrals.
Nonetheless, even if the difficulty of generating such a high-order mesh is lower than building fully analysis-suitable boundary-conforming parameterizations, it still remains a challenging and time consuming task for complex geometries.

An appealing alternative to these two approaches is the use of \emph{moment fitting} techniques~\cite{joulaian2016numerical,Hubrich_2017,Hubrich_2019,bui2020efficient} in which coarse, but accurate, quadrature rules are generated for complex integration domains by tuning the positions and/or weights of the quadrature points.
Nevertheless, these methods come at a price: The creation of tailored quadrature rules requires the computation of polynomial integrals over complex domains at a pre-processing stage, what calls for the use of alternative integration techniques.

Finally, there exists a fourth group of strategies for computing integrals over curved polyhedrons that lies in deriving dedicated integration rules for specific classes of integrands, as for instance polynomial functions.
Indeed, it is known that integrating polynomials and other homogeneous functions over (curved) polyhedrons can be done more efficiently by invoking the divergence theorem~\cite{Lasserre_1998,Gonzalez_Ochoa_1998,Mousavi_2010,Chin_2015,Chin_2020}.
These results can be exploited in several ways: One can perform a \emph{polynomial approximation} of the integrands of interest such that the integration can be done straightforwardly~\cite{Ventura_2006,Duczek_2015,Abedian_2019}; those specific rules can be applied at the pre-processing stage of \emph{moment-fitting} methods~\cite{Mueller_2013,Hubrich_2017,Hubrich_2019}; or by invoking other specific procedures~\cite{Sudhakar_2014,Gunderman_2021}.

Within this category, worth mentioning are the recent works~\cite{Chin_2020,Gunderman_2021}, where the divergence theorem is used for transforming volumetric integrals into either surface or line integrals.
In~\cite{Gunderman_2021}, the authors reduced 3D integrals of general functions to 1D integrals, that are finally evaluated using fine quadrature rules.
Similarly, in~\cite{Chin_2020} the complexity of 3D integrals is reduced to just vertices evaluations in the case of planar polyhedra.
While for the case of B-reps composed of B\'ezier triangles or non-trimmed B-splines patches, the authors in~\cite{Chin_2020} applied the divergence theorem just once, transforming 3D integrals in 2D ones, that are approximated through standard quadrature rules.

Aligned with these ideas, in this work we present a fully \emph{quadrature-free} method for integrating polynomials over general B-rep models enclosed by trimmed spline surfaces.
The procedure is based on two successive applications of the divergence theorem, reducing volumetric integrals to first surface and then line integrals, that are computed analytically up to machine precision.
Hence, this can be seen as a generalization of those precedent works, eliminating the need of quadrature rules.
Furthermore, we show how this integration procedure, combined with a consistent polynomial approximation step, leads to a new analysis tool for immersed isogeometric methods that skips the need of complex quadrature rules.


The developed approach is presented as follows: We firstly introduce in Section~\ref{sec:immersed_iga} the basics regarding immersed isogeometric analysis to further detail the scope of application of this work, and describe a consistent approximation step required for transforming the involved integrands into polynomials.
Then, in Section~\ref{sec:BRep}, we discuss the geometric modeling via splines, trimming, and boundary-representation, as commonly undertaken in CAD.
In Section~\ref{sec:integration}, the proposed \emph{quadrature-free} integration over B-Reps is presented.
Finally, in Section~\ref{sec:numericalExamples}, we solve elliptic PDEs and perform several numerical experiments to confirm the accuracy of the approach.
Lastly, concluding remarks are summarized in Section~\ref{sec:conclusions}.

\begin{figure*}[tp]
  \centering
  \includegraphics[width=0.75\linewidth]{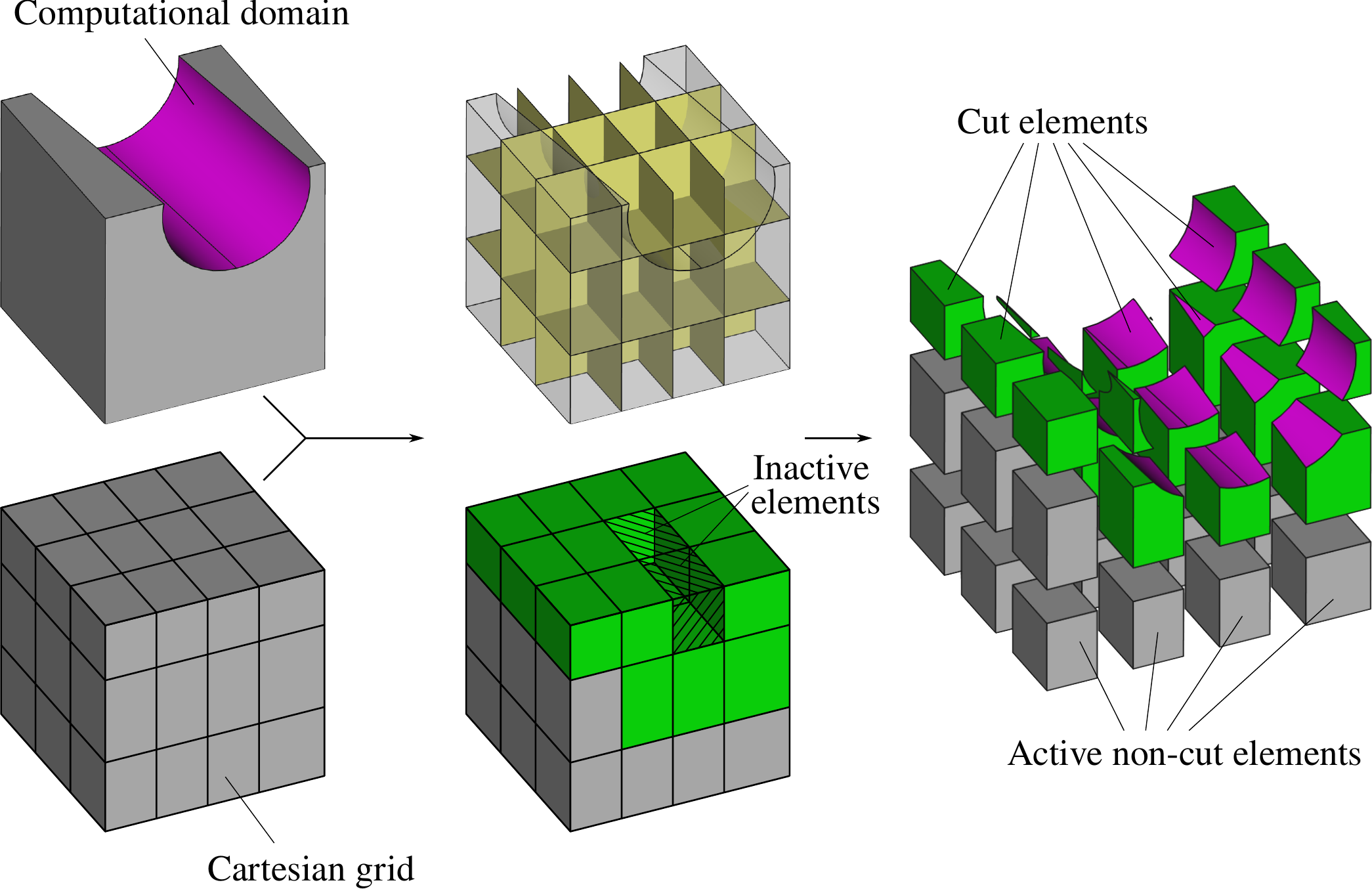}
  \caption{Immersed methods setting.}
  \label{fig:immersedIGA}
\end{figure*}

\section{Immersed isogeometric analysis} \label{sec:immersed_iga}

With the aim of introducing immersed methods, the used notation, and the main ideas behind this work, let us first introduce a classical Poisson's problem as our driving example.
Even if the problem is presented in a 3D context, the same ideas are directly applicable to 2D problems.

Let $\Omega\subset\mathbb{R}^3$ be the computational domain whose boundary is partitioned as $\Gamma_N\cup\Gamma_D=\partial\Omega$ and $\Gamma_N\cap\Gamma_D=\emptyset$.
We also define a functional space $H^1_D(\Omega)=\{ v\in H^1(\Omega) \ :\  v\vert_{\Gamma_D}  =0 \}$, such that the Poisson's problem reads: find $u\in{H^1_D(\Omega)}$ solution of:
\begin{equation}
  \begin{split}
    -\nabla\cdot(\diffusivityCoef\nabla{u}) = f & ~~\text{in}~\Omega, \\
    \nabla{u}\cdot\bn=g & ~~\text{on}~\Gamma_N, \\
    u=0 & ~~\text{on}~\Gamma_D,
  \end{split}
  \label{eq:strongpoisson}
\end{equation}
where $\diffusivityCoef\in L^2(\Omega)^{3\times{3}}$ is the symmetric diffusivity operator; $f\in L^2(\Omega)$ and $g\in H^{-1/2}(\Gamma_N)$ are the source and Neumann terms, respectively; and~$\bn\in\mathbb{R}^3(\partial\Omega)$ is the outward pointing unit normal on the boundary.
For the sake of clarity, and without constituting any limitation, in the problem~\eqref{eq:strongpoisson} and hereinafter we assume homogeneous Dirichlet boundary conditions.

The associated weak problem can be written as: find $u\in H^1_D(\Omega)$ such that
\begin{equation}
  a(u,v) = {b}(v), \quad\forall{v}\in H^1_D(\Omega)\,,
\end{equation}
where
\begin{equation} \label{eq:poissonweak}
  \begin{split}
    a(u,v)
    &=
    \int_{\Omega} \nabla{u}\cdot\diffusivityCoef\nabla{v} \,\diff\Omega,\\
    b(v)
    &=
    \int_{\Omega} f\,{v} \,\diff\Omega
    +
    \int_{\Gamma_N} g\,{v} \,\diff\Gamma.
  \end{split}
\end{equation}

\subsection{Immersed methods} \label{sec:immersed}


The philosophy behind immersed methods is depicted in Figure~\ref{fig:immersedIGA}.
It consists in embedding the computational domain~$\Omega$ into a grid $\T_h(\Omega_0)$ of a larger domain $\Omega_0$, such that $\Omega\subset\Omega_0\subset\mathbb{R}^3$.
The solution of the weak problem~\eqref{eq:poissonweak} is then discretized over a subset of the grid $\T_h(\Omega_0)$, what allows to decouple the solution discretization from the actual geometry.
This simple and rather straightforward procedure is the one and only mesh generation task to undertake within immersed-like approaches, making this class of methods very appealing.
Indeed, this can largely ease the design-to-analysis workflow since the computational domain can be directly prescribed as a geometric model with any representation commonly used in CAD, as for instance the Boundary-Representation (detailed in Section~\ref{sec:BRep}).
In return, the price to pay during the analysis lies in the introduction of so-called cut or trimmed elements, as illustrated in Figure~\ref{fig:immersedIGA}.
The presence of these elements raises a major difficulty: the integration of quantities over cut elements (as discussed in the introduction, see Section~\ref{sec:intro}).
This work focuses in this particular challenge one would face when dealing with enriched or unfitted finite element methods over B-Rep models.

As the computational domain is $\Omega$ and not $\Omega_0$, the partition $\T_h(\Omega_0)$ is restricted to a subset $\T_h(\Omega)$ as:
\begin{equation}
  \T_h(\Omega)  : =   \{ Q\ \vert\ \forall Q\in \T_h({\Omega_0})  \ :\ Q\cap \Omega \neq \emptyset    \}.
\end{equation}
Indeed, the grid $\T_h(\Omega_0)$ naturally splits the domain $\Omega_0$ into three complementary partitions of elements:
\begin{subequations}
\begin{align}
  \T^\Gamma_h(\Omega)  : &=   \{ Q\ \vert\ \forall Q\in \T_h({\Omega})  \ :\ Q\cap \Omega \neq Q  \}\,,\\
  \T^{\text{int}}_h(\Omega)  : &=   \{ Q\ \vert\ \forall Q\in \T_h({\Omega})  \ :\ Q\cap \Omega = Q    \}\,,\\
  \T^0_h(\Omega_0)  : &=   \{ Q\ \vert\ \forall Q\in \T_h({\Omega_0})  \ :\ Q\cap \Omega = \emptyset    \}\,,
\end{align}
\end{subequations}
such that $\T_h(\Omega)=\T^{\text{int}}_h(\Omega)\cup\T^\Gamma_h(\Omega)$ and $\T_h(\Omega_0)=\T_h(\Omega)\cup\T^0(\Omega_0)$.
As depicted in Figure~\ref{fig:immersedIGA}, the elements belonging to these three subsets are denoted as cut, non-cut, and inactive elements, respectively.

In this work we limit our discussion to the case of 3D immersed isogeometric methods, nevertheless, the presentation is kept rather general and can be easily adapted to generic immersed methods~\cite{peskin2002immersed} or particular cases as, for instance, CutFEM~\cite{burman2015cutfem} or Finite Cell Methods~\cite{parvizian2007finite}, among others.

In order to solve numerically the weak problem~\eqref{eq:poissonweak} we construct a discrete spline space $\mathbb{V}_h(\Omega_0)$ over the grid $\T_h(\Omega_0)$ as:
\begin{equation}
  \mathbb{V}_h(\Omega_0)  =   \text{span}\{ N^p_i,~i\in\mathcal{I}_0\}\,,
\end{equation}
where $N^p_i$ denotes generic spline basis functions of degree $p>0$ and arbitrary continuity (up to $p-1$),  and $\mathcal{I}_0$ is the set of indices of those basis functions, such that $\text{dim}(\mathbb{V}_h(\Omega_0)) = \# \mathcal{I}_0$.
In this work we use tensor-product B-splines, but the extensions to other cases as, \emph{e.g.}, hierarchical splines~\cite{giannelli2012thb} or T-splines~\cite{bazilevs2010isogeometric}, is straightforward.
For the sake of simplicity, henceforward we drop the superscript $p$ from $N^p_i$ and assume that the spline degree $p$ is constant along the three parametric directions.

The support of some basis functions of the space $\mathbb{V}_h(\Omega_0)$ may not intersect the domain $\Omega$ and, consequently, they do not contribute to the solution of the problem~\eqref{eq:poissonweak}.
Therefore, we trim the space $\mathbb{V}_h(\Omega_0)$ as:
\begin{equation} \label{eq:space}
  \mathbb{V}_h(\Omega)  =   \text{span}\{ N_i\in \mathbb{V}_h(\Omega_0)\ :\  \text{supp}\{ N_i \} \cap \Omega \neq \emptyset\}\,,
\end{equation}
that, as already studied in~\cite{Antolin_2019}, holds optimal approximation properties.
It is a well-known fact that the active support of some basis functions in~$\mathbb{V}_h(\Omega)$ ($\text{supp}\{ N_i \} \cap \Omega$) may be small, what could yield ill-conditioned operators.
This is an active research topic~\cite{B_chet_2005,de_Prenter_2017,Marussig_2017,Buffa_2020} that exceeds the scope of this work.

Henceforward, we assume the Dirichlet boundary $\Gamma_D$ to be such that $\Gamma_D\subset\partial\Omega_0\cap\partial\Omega$, what grants the strong enforcement of Dirichlet boundary conditions.
The opposite case ($\Gamma_D\not\subset\partial\Omega_0$) entails the imposition of Dirichlet conditions in a weak sense.
We refer the interested reader to~\cite{Hansbo_2002,Ruess_2013,Pande_2021} for a dedicated discussion and to~\cite{Buffa_2020} for a study,  in the case of spline spaces, of the inherent stability issues.

Thus, by means of the assumption $\Gamma_D\subset\partial\Omega_0\cap\partial\Omega$, we can define the space:
\begin{equation} \label{eq:spaceD}
  \mathbb{V}^D_h(\Omega)  =   \{  v_h\in \mathbb{V}_h(\Omega)  \ :\ v_h\,\vert_{\Gamma_D} =0 \}\,.
\end{equation}
that allows us to discretize the continuous weak problem~\eqref{eq:poissonweak} as: find~${u}_h\in\mathbb{V}^D_h(\Omega)$ solution of:
\begin{equation} \label{eq:weakdiscrete}
  a({u}_h,{v}_h) = {b}({v}_h), \quad\forall{v}_h\in{\mathbb{V}^D_h(\Omega)},
\end{equation}
where the discrete versions of the bilinear form~${a}$ and the linear form~${b}$ are decomposed as:
\begin{equation} \label{eq:discreteforms}
  \begin{split}
    a(u_h,v_h) =&
    \sum_{Q\in\T^{\text{int}}_h(\Omega)} \int_{Q} \nabla{u_h}\cdot\diffusivityCoef\nabla{v_h} \,\diff Q \\
    &+
    \sum_{Q\in\T^\Gamma_h(\Omega)} \int_{Q\cap\Omega} \nabla{u_h}\cdot\diffusivityCoef\nabla{v_h} \,\diff Q\,, \\
    b(v_h) =&
    \sum_{Q\in\T^{\text{int}}_h(\Omega)} \int_{Q} f\,{v_h} \,\diff Q \\
    &+
    \sum_{Q\in\T^{\Gamma}_h(\Omega)} \int_{Q\cap\Omega} f\,{v_h} \,\diff Q \\
    &+
    \sum_{Q\in\T^{\Gamma}_h(\Omega)} \int_{Q\cap\Gamma_D} g\,{v_h} \,\diff \Gamma\,.
  \end{split}
\end{equation}

The computation of the integrals over non-cut elements $Q\in\T^{\text{int}}_h(\Omega)$ is straightforward and can be performed using classical quadrature schemes.
However, the evaluation of integrals over cut elements $Q\in\T^{\Gamma}_h(\Omega)$ is a challenging problem and one of the Achilles' heels of isogeometric immersed methods in 3D (see the related discussion in Section~\ref{sec:intro}).
The main contribution of this article regards the computation of those integrals through a quadrature-free approach for the case of cut elements defined as B-Rep models.
This procedure is presented in Section~\ref{sec:integration}.
Nonetheless, this method is only applicable to the case in which the integrands are polynomial functions.
Thus, before introducing it, in the next section the integrals in~\eqref{eq:discreteforms} are transformed such as they only rely on polynomial integrands.

\subsection{Polynomial approximation of finite element operators} \label{sec:polapprox}

When considering spline discretizations over the grid~$\T_h(\Omega)$, the terms~$\nabla{u_h}$, $\nabla{v_h}$, and $v_h$ in the operators~\eqref{eq:discreteforms} take polynomial forms $\forall Q\in\T_h(\Omega)$.
On the contrary, the datum quantities involved (\emph{i.e.}, $\diffusivityCoef$, $f$, and~$g$) may not be polynomials in general.

Hence, in order to work with integrals that only present polynomial integrands, we seek to exploit a key result introduced in~\cite{Mantzaflaris_2015}: It is possible to perform a polynomial approximation of the integrands in~\eqref{eq:discreteforms} without deteriorating the solution.
More specifically, instead of solving the problem~\eqref{eq:weakdiscrete}, we consider the following approximate problem: find~$\bar{u}_h\in\mathbb{V}^D_h(\Omega)$ solution of:
\begin{equation}
  \bar{a}(\bar{u}_h,{v}_h) = \bar{b}({v}_h), \quad\forall{v}_h\in\mathbb{V}^D_h(\Omega),
  \label{eq:weakpbapprox}
\end{equation}
where the discrete forms in~\eqref{eq:discreteforms} are replaced by:
\begin{equation} \label{eq:weakpbapprox_forms}
  \begin{split}
    \bar{a}(\bar{u}_h,v_h) &=
    \sum_{Q\in\T^{\text{int}}_h(\Omega)} \int_{Q} \nabla\bar{u}_h\cdot\bar{\diffusivityCoef}\nabla{v_h} \,\diff Q \\
    &+
    \sum_{Q\in\T^\Gamma_h(\Omega)} \int_{Q\cap\Omega} \nabla\bar{u}_h\cdot\bar{\diffusivityCoef}\nabla{v_h} \,\diff Q\,,\\
    \bar{b}(v_h) &=
    \sum_{Q\in\T^{\text{int}}_h(\Omega)} \int_{Q} \bar{f}\,{v_h} \,\diff Q \\
    &+
    \sum_{Q\in\T^{\Gamma}_h(\Omega)} \int_{Q\cap\Omega} \bar{f}\,{v_h} \,\diff Q \\
    &+
    \sum_{Q\in\T^{\Gamma}_h(\Omega)} \int_{Q\cap\Gamma_D} \bar{g}\,{v_h} \,\diff \Gamma\,,
  \end{split}
\end{equation}
that involves the following polynomial approximations:
\begin{equation}
  \bar{\diffusivityCoef} = \Pi^{h}{\diffusivityCoef},\qquad
  \bar{f} = \Pi^{h}{f},\qquad
  \bar{g} = \Pi^{h}{g}.
  \label{eq:polynomialapprox}
\end{equation}

In the approximations above, the projection spaces must be chosen carefully, such that the introduced consistency errors do not pollute the numerical solution.
Thus, by recalling \cite[Theorem~13]{Mantzaflaris_2015}, we know that the projection of $\diffusivityCoef$, $f$, and~$g$ into spline spaces of degree~$q\geq{p-1}$ yields a solution~$\bar{u}_h$ that approximates optimally the true solution $u$, presenting convergence order $p$ for the error measured in the $H^1$~semi-norm when the mesh size $h\to 0$.
In \cite{Mantzaflaris_2015}, the authors also observed, through numerical experiments, that a projection degree $q>p-1$ yields optimal convergence order also respect to the $L^2$ norm of the error (rate $p+1$).
\begin{myremark} \label{rmk:mapping}
  The non-polynomial nature of the quantities $\diffusivityCoef$, $f$, and~$g$ may derive from an additional mapping that further deforms the domain $\Omega_0$ (see, \emph{e.g.}, \cite{Antolin_2019b}).
  A numerical example addressing this case is presented in Section~\ref{sec:poisson2D} (the multi-perforated quarter of annulus).
  On the contrary, these quantities might be low-order polynomials (even zero-order polynomials) by construction and it is therefore not necessary to project them into polynomial spaces.
\end{myremark}

In~\cite{Mantzaflaris_2015}, the projections~\eqref{eq:polynomialapprox} are performed patch-wise.
Nevertheless, the same error estimates hold in the case they are carried out in an element-wise way, that is the case of this work.
This results in polynomial approximations that are element-wise discontinuous.
Thus, for each element $Q\in\T_h(\Omega)$ we introduce a local $L^2$-projector:
\begin{equation}
  \Pi^{h}_{Q}:{L}^2(Q)
  \to \mathbb{Q}_{q,\,q,\,q}(Q)\,,\quad \forall Q\in\T_h(\Omega)\,,
  \label{eq:projectorDef}
\end{equation}
where $\mathbb{Q}_{q_1,q_2,\dots,q_{m}}$ denotes the space of tensor-product polynomials with degrees $(q_1,q_2,\dots,q_{m})$ along the $m$ parametric directions.


By employing a tensor-product Bernstein basis, the projected quantities $\bar{\diffusivityCoef}$, $\bar{f}$, and~$\bar{g}$ restricted to element $Q$ can be expressed as:
\begin{equation}
  \begin{split}
    \left.\bar{\diffusivityCoef}\right\rvert_{Q} = \sum^{\left(q+1\right)^3}_{k=1} B^{\mathbf{q}}_k\,\bar{K}^{(Q)}_k,\\
    \left.\bar{f}\right\rvert_{Q} = \sum^{\left(q+1\right)^3}_{k=1} B^{\mathbf{q}}_k\,\bar{f}^{(Q)}_k,\\
    \left.\bar{g}\right\rvert_{Q} = \sum^{\left(q+1\right)^3}_{k=1} B^{\mathbf{q}}_k\,\bar{g}^{(Q)}_k\,,
  \end{split}
  \label{eq:polynomialapprox2}
\end{equation}
where $\bar{K}^{(Q)}_k\in\mathbb{R}^{3\times{3}}$, $\bar{f}^{(Q)}_k\in\mathbb{R}$, and~$\bar{g}^{(Q)}_k\in\mathbb{R}$ are the projection coefficients, and $B^{\mathbf{q}}_k$ are tensor-product Bernstein polynomials defined over $Q$ and with degrees $\mathbf{q}=(q,q,q)$ such that
\begin{equation} \label{eq:polspaces}
  \mathbb{Q}_{q,\,q,\,q}(Q) = \text{span}\{B^{\mathbf{q}}_{k}\ :\  k=1,\dots,\left(q + 1\right)^3\}\,.
\end{equation}
We refer the interested reader to Appendix~\ref{sec:appA} for a discussion about Bernstein polynomials.

\begin{figure*}[t]
  \centering
  \includegraphics[width=0.75\linewidth]{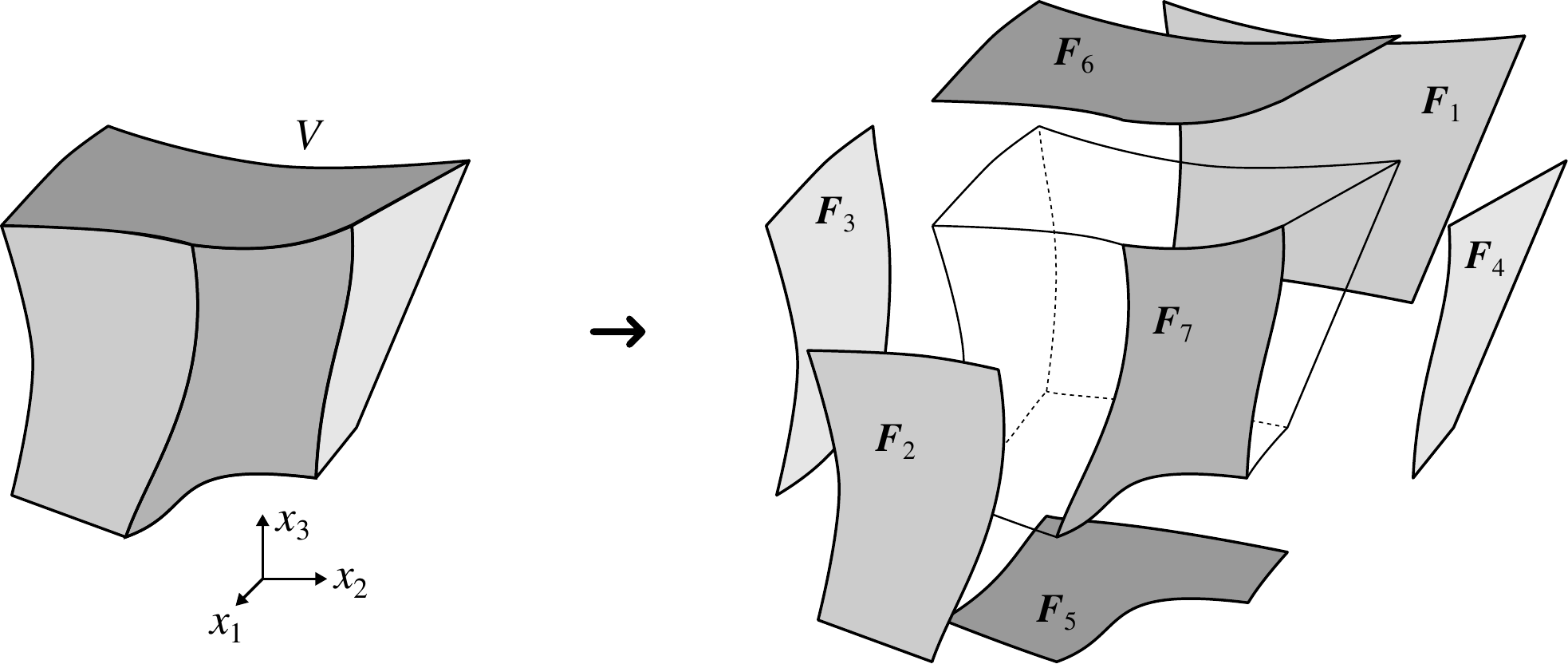}
  \caption{Boundary representation of a volumetric domain~$V$.}
  \label{fig:brepVolbySurf}
\end{figure*}

\subsection{Operators assembly through lookup tables}
In what follows, we detail the assembly of the elemental stiffness matrix and the right-hand-side vector associated to the operators~\eqref{eq:weakpbapprox_forms}.
Thus, plugging the projections~\eqref{eq:polynomialapprox2} into~\eqref{eq:weakpbapprox_forms}, a single entry of the elemental matrix and vector can be computed as:
\begin{equation} \label{eq:localoperators}
  \begin{split}
    \mathsf{A}_{ij}^{(Q)}
    &= \sum_{k=1}^{(q+1)^3} \bar{K}^{(Q)}_k : \int_{Q\cap\Omega} B^{\mathbf{q}}_k \big(\nabla{N_i} \otimes \nabla{N_j}\big) \,\diff Q,
    \\
    \mathsf{b}_{i}^{(Q)}
    &= \sum_{k=1}^{(q+1)^3} \bar{f}^{(Q)}_k \int_{Q\cap\Omega} B^{\mathbf{q}}_k N_i \,\diff Q \\
    &+ \sum_{k=1}^{(q+1)^3} \bar{g}^{(Q)}_k \int_{Q\cap\Gamma_D} B^{\mathbf{q}}_k N_i \,\diff \Gamma,
  \end{split}
\end{equation}
where~$N_i,N_j\in\mathbb{V}(\Omega)$ are test and trial basis functions, respectively.
In the expressions above it is easy to realize that all the integrands restricted to a single element $Q$ are polynomials:
\begin{subequations} \label{eq:termsdegrees}
  \begin{align}
    &\left.B^{\mathbf{q}}_k \big(\nabla{N_i}\otimes\nabla{N_j}\big)\right\rvert_{Q} \in \mathbb{Q}_{2 p +q,\,2 p +q,\,2 p +q}(Q)\,,\\
    &\left.B^{\mathbf{q}}_k N_i\right\rvert_{Q} \in \mathbb{Q}_{p +q,\,p +q,\,p +q}(Q)\,.
\end{align}
\end{subequations}
Notice also that the functions~$N_i$, $N_j$, and~$B^{\mathbf{q}}_k$ are naturally defined over the full support of each element $Q$, and not only over its active part $Q\cap\Omega$.

Finally, by exploiting their polynomial nature, the element integrals in~\eqref{eq:localoperators} can be computed as:
\begin{subequations}
  \begin{align}
    \int_{Q\cap\Omega} B^{\mathbf{q}}_k\big(\nabla{N_i}&\otimes\nabla{N_j}\big) \,\diff Q \nonumber \\
                                                       &= 
                                                         \sum^{\left(2 p + q + 1\right)^3}_{\alpha=1}  \mathbf{\mathsf{K}}^{(Q)}_{i,j,k,\alpha} \int_{Q\cap\Omega} B^{\mathbf{r}}_{\alpha} \diff Q\\
    \label{eq:integral_source}
    \int_{Q\cap\Omega} B^{\mathbf{q}}_k N_i \,\diff Q
    &=
      \sum^{\left(p + q + 1\right)^3}_{\beta=1} \mathsf{F}^{(Q)}_{i,k,\beta} \int_{Q\cap\Omega} B^{\mathbf{s}}_{\beta} \diff Q\\
    \int_{Q\cap\Gamma_D} B^{\mathbf{q}}_k N_i \,\diff Q
                                                       &=
                                                         \sum^{\left(p + q + 1\right)^3}_{\beta=1} \mathsf{G}^{(Q)}_{i,k,\beta} \int_{Q\cap\Gamma_D} B^{\mathbf{s}}_{\beta} \diff \Gamma
  \end{align}
\end{subequations} 
where $B^{\mathbf{r}}_{\alpha}$ and $B^{\mathbf{s}}_{\beta}$ are tensor-product Bernstein polynomials with degrees $\mathbf{r}=(2 p + q,\,2 p + q,\,2 p + q)$ and $\mathbf{s}=(p + q,\,p + q,\,p + q)$.
$\mathbf{\mathsf{K}}^{(Q)}_{i,j,k,\alpha}\in\mathbb{R}^{3\times 3}$ and $\mathsf{F}^{(Q)}_{i,k,\beta},\,\mathsf{G}^{(Q)}_{i,k,\beta}\in\mathbb{R}$ are element dependent constant coefficients that can be calculated by means of the B\'ezier extraction operators~\cite{borden2011isogeometric,d2018multi,scott2011isogeometric} associated to the spline space~$\mathbb{V}_h(\Omega)$.

Then, the assembly of the operators~\eqref{eq:localoperators} reduces to the computation of the coefficients $\mathbf{\mathsf{K}}^{(Q)}_{i,j,k,\alpha}$, $\mathsf{F}^{(Q)}_{i,k,\beta}$, and $\mathsf{G}^{(Q)}_{i,k,\beta}$, as well as the integrals\footnote{Due to the fact that $\mathbb{Q}_{p +q,\,p +q,\,p +q}\subset\mathbb{Q}_{2 p +q,\,2 p +q,\,2 p +q}$, the integrals $\int_{Q\cap\Omega} B^{\mathbf{s}}_{\beta} \diff Q$ in~\eqref{eq:integral_source} can be computed as linear combinations of the integrals $\mathsf{I}^{3\textup{D}}_{Q,\alpha}.$}:
\begin{equation} \label{eq:integrals}
  \mathsf{I}^{3\textup{D}}_{Q,\alpha} = \int_{Q\cap\Omega} B^{\mathbf{r}}_{\alpha} \diff Q\,,\quad
  \mathsf{I}^{2\textup{D}}_{Q,\beta} = \int_{Q\cap\Gamma_D} B^{\mathbf{s}}_{\beta} \diff \Gamma\,.
\end{equation} 
Thus, the integrals $\mathsf{I}^{3\textup{D}}_{Q,\alpha}$ and $\mathsf{I}^{2\textup{D}}_{Q,\beta}$ can be precomputed for every element $Q$ and stored in lookup tables, that will be accessed along the assembly process to create the elemental operators, in a similar way as proposed in~\cite{Mantzaflaris_2015}.

Nevertheless, as discussed in Section~\ref{sec:intro}, the computation of the integrals~\eqref{eq:integrals} is a challenging task.
In the case of non-cut elements, their evaluation is straightforward: It can be precomputed analytically for a single unit cube and subsequently adapted to every non-cut element's domain through simple transformations (translations and scalings).
But in the case of cut elements the evaluation of the integrals $\mathsf{I}^{3\textup{D}}_{Q,\alpha}$ and $\mathsf{I}^{2\textup{D}}_{Q,\beta}$ is far from simple.
For that purpose, in Section~\ref{sec:integration} we propose a quadrature-free approach for the common case in which the active part of elements ($Q\cap\Omega$) can be defined through a B-Rep, discussed in Section~\ref{sec:BRep}.

\section{Geometric modeling via boundary representation}
\label{sec:BRep}

In this section we introduce the notation and some basic concepts about splines and geometric modeling.
Hence, we provide a mathematical way of describing the active part of the cut elements $Q\cap\Omega$, discussed in the previous section,  by means of B-Rep representations.
This constitutes the basis for the integration method presented in Section~\ref{sec:integration}.

\begin{figure*}[t]
  \centering
  \includegraphics[width=0.65\linewidth]{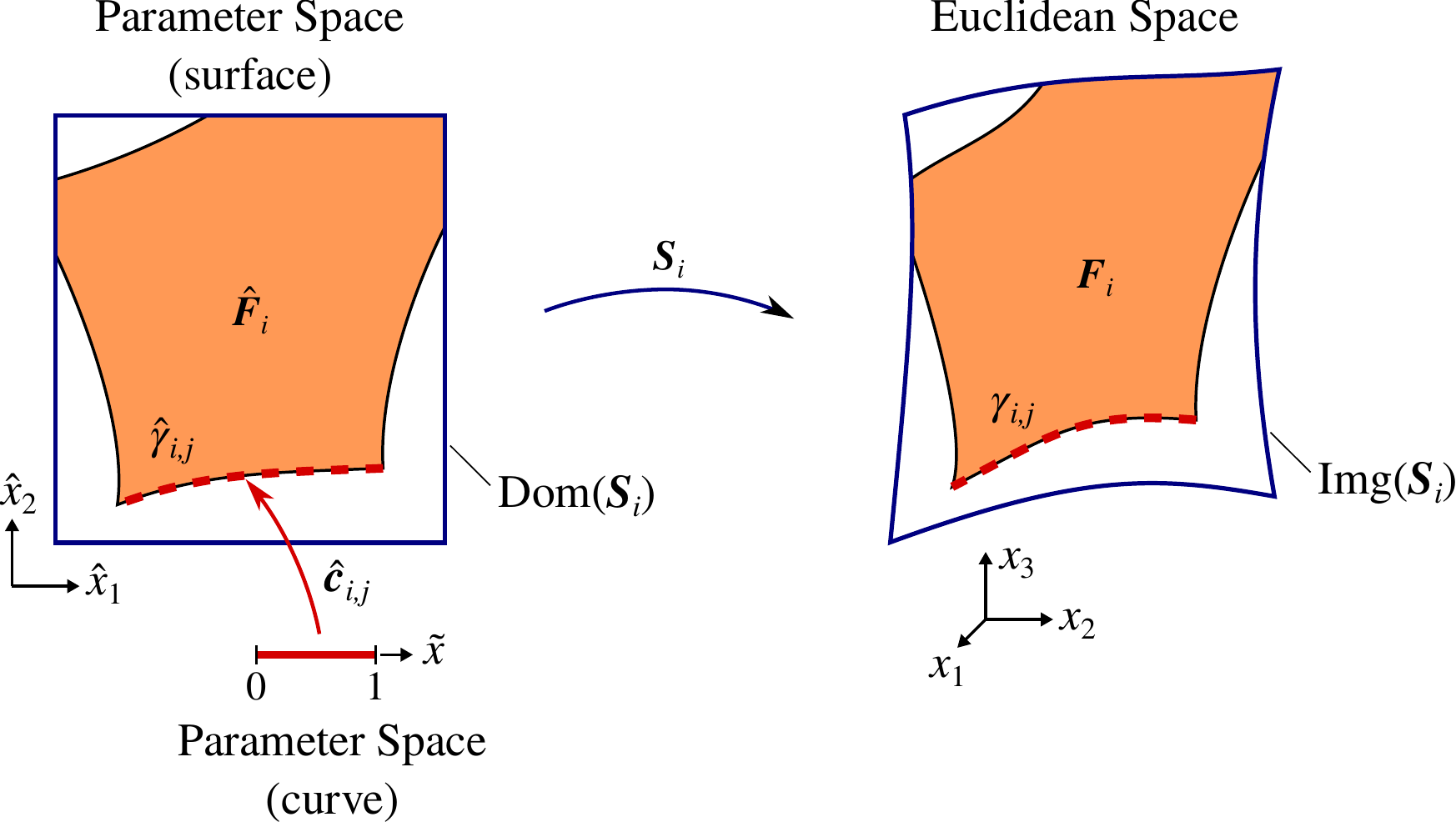}
  \caption{Description of the involved geometrical entities in the definition of trimmed parametric surfaces.}
  \label{fig:fig_2}
\end{figure*}

\subsection{Spline representation}

Splines are considered a \emph{de facto} standard in Computer-Aided Design and have been extensively studied in the literature, see for instance~\cite{Cohen2001,Farin2001,Piegl_1997}.
Among the different representation techniques available, in this work we focus on the use of polynomial mappings, and more specifically, B-spline and B\'ezier curves and surfaces.
A B-spline or B\'ezier curve $\bm{c}$ can be expressed in the form:
\begin{equation} \label{eq:splineCurveMap}
  \bm{c}:[0,1] \to \mathbb{R}^d,~~
  \tilde{x}\mapsto \boldsymbol{c}(\tilde{x}) = \sum_{i=1}^{n} N^{p}_i(\tilde{x}) \boldsymbol{P}_{i}\,,
\end{equation}
where $N^{p}_i$ are univariate basis functions, either B-splines or Bernstein polynomials, of degree $p$, and $\boldsymbol{P}_i\in\mathbb{R}^d$ are their associated control points, being $d$ the space dimension.
In Appendix~\ref{sec:appA} we provide further details about Bernstein polynomials and B\'ezier geometries, that are extensively used in this work.
For an in-depth discussion about B-Splines, we refer the interested reader to the existing literature \cite{Cohen2001,Farin2001,Piegl_1997}.

Using tensor-product combinations of those basis functions, B-Spline and B\'ezier surfaces $\bm{S}$ can be constructed as:
\begin{equation} \label{eq:splineSurfaceMap}
  \begin{split}
    \bm{S}:~~[0,1]^{2} &\to \mathbb{R}^d,\\
    (\hat{x}_1,\hat{x}_2)&\mapsto \sum_{i=1}^{n_1}\sum_{j=1}^{n_2} N^{p_1}_i(\hat{x}_1)N^{p_2}_j(\hat{x}_2) \boldsymbol{P}_{i,j},
  \end{split}
\end{equation}
where $N^{p_1}_i$ and $N^{p_2}_j$ are univariate B-spline or Bernstein basis functions of degrees $p_1$ and $p_2$, respectively, and $\boldsymbol{P}_{i,j}\in\mathbb{R}^d$ are the associated control points.
For the sake of simplicity, we assumed the parametric domains of the mappings~\eqref{eq:splineCurveMap} and \eqref{eq:splineSurfaceMap}, $\operatorname{Dom}(\bm{c})$ and $\operatorname{Dom}(\bm{S})$, to be $[0,1]$ and $[0,1]^2$, respectively.

\subsection{Trimmed surfaces and boundary representations} \label{ssec:brep}

Simple spline mappings~\eqref{eq:splineCurveMap} and \eqref{eq:splineSurfaceMap} cannot represent complex real-world geometries.
Instead, multitude of these geometric objects are usually combined for such a purpose.
More specifically, Boolean operations (namely, unions, differences, and/or intersections) of several geometrical entities are commonly adopted in Computer-Aided Design~\cite{Cohen2001}.
By means of these operations, volumetric geometries are often represented in an implicit way: the volume enclosed by a set of, possibly trimmed, boundaries surfaces.
This paradigm, known as Boundary Representation (B-Rep)~\cite{requicha1992solid,braid1973designing} and extensively used in industrial modeling tools, is considered throughout this work.

As illustrated in Figure~\ref{fig:brepVolbySurf}, we consider a domain~$V\subset\mathbb{R}^3$, non-simply connected in general, whose boundary~$\partial V$ is defined by a set of connected faces $\Fi,~i=1,\dots,n_F$, such as:
\begin{equation} \label{eq:brepPhysical}
  \partial V=\cup_{i=1}^{n_F}\Fi.
\end{equation}
The domain $V$ may correspond to the active part of the cut elements $Q\cap\Omega$ discussed in Section~\ref{sec:immersed}.

\begin{figure*}[t]
  \centering
  \includegraphics[width=0.50\linewidth]{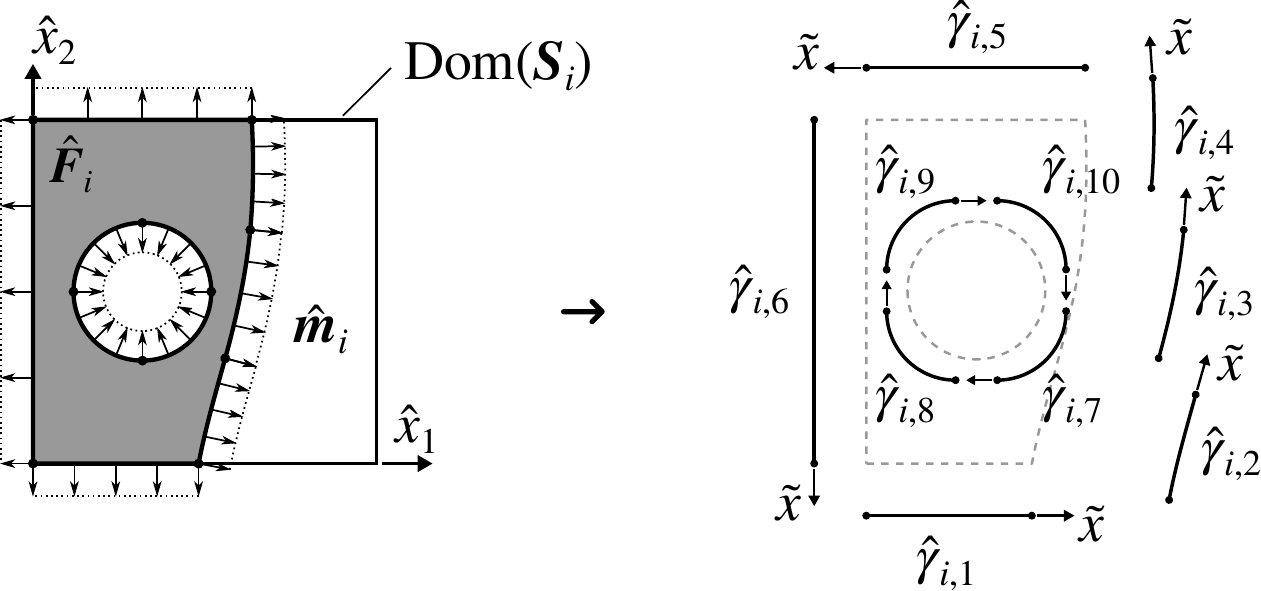}
  \caption{Boundary representation of trimmed faces.
           External boundaries follow a counter-clockwise orientation while the internal ones are clockwise oriented.}
  \label{fig:surfaceTOcurves}
\end{figure*}

We consider the faces~$\Fi$ to be defined as trimmed B-spline or B\'ezier surfaces that are piecewise smooth.
Every trimmed face~$\Fi$ is composed of two elements: an underlying spline surface mapping~$\Si$ of the form~\eqref{eq:splineSurfaceMap}, and a group of connected curvilinear segments $\ghij\subset\operatorname{Dom}(\Si),~j=1,\dots,n_{c,i}$, that delimit the active region of $\operatorname{Dom}(\Si)$ (see Figures~\ref{fig:fig_2} and \ref{fig:surfaceTOcurves}).
We denote this active region as $\Fhi\subset\operatorname{Dom}(\Si)$.

Each segment $\ghij$ is the image of a spline curve mapping~$\chij:[0,1] \to \ghij$ of the form~\eqref{eq:splineCurveMap}.
Thus, the boundary of the active region $\Fhi$ is:
\begin{equation}
  \label{eq:Fhi}
  \begin{split}
    \partial\Fhi&= \cup_{j=1}^{n_{c,i}} \ghij,\\
    \ghij&=\{ \bxh\in\mathbb{R}^2\ \vert\ \xt\in[0,1]\ :\ \bxh=\chij(\xt) \},
  \end{split}
\end{equation}
therefore, we can define~$\Fi$ as:
\begin{equation}
   \Fi=\{ \bx\in\mathbb{R}^3 \ \vert\ \bxh\in\Fhi \ :\ \bx=\Si(\bxh) \}.
\end{equation}
We again refer to Figure~\ref{fig:fig_2} where all the introduced quantities are depicted for an illustrative example.

\begin{figure*}[t]
  \centering
  \includegraphics[width=0.65\linewidth]{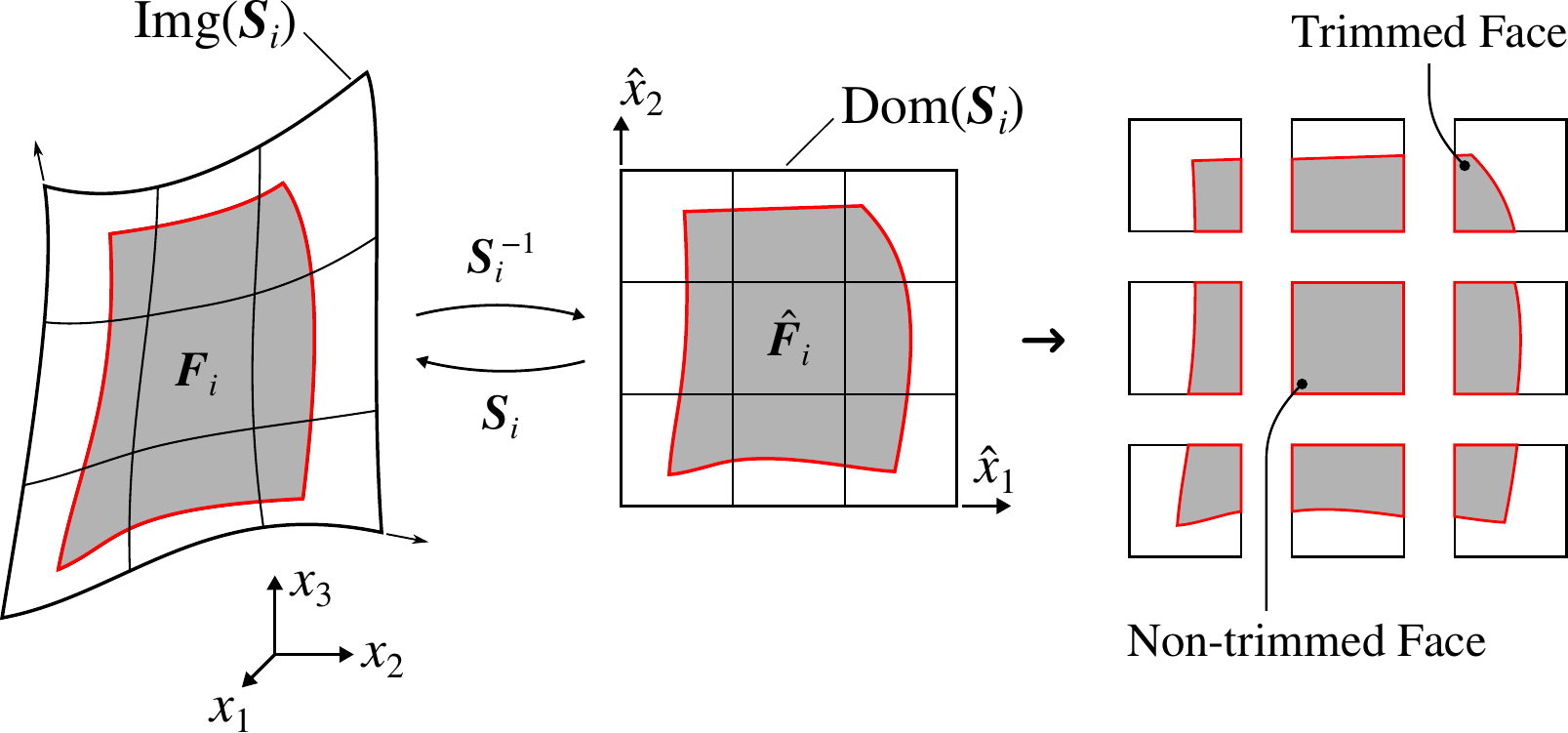}
  \caption{Split of a trimmed B-Spline surface into trimmed B\'ezier surfaces.}
  \label{fig:splitPiecewisePolySurf}
\end{figure*}

\begin{myremark} \label{rmk:rationals}
In order to work exclusively with pure polynomial representations, instead of (rational) piecewise polynomials, in this work we only consider non-rational B\'ezier curves and surfaces.
Using only B\'eziers does not constitute any limitation: By refining at its internal knots, any face $\Fi$, defined by means of B-spline curves and surfaces, can be easily split into a set of trimmed B\'ezier faces, whose underlying curves and surfaces are B\'eziers (see Figure~\ref{fig:splitPiecewisePolySurf}).
On the other hand, the exclusive use of non-rational polynomials may be a limiting factor as it turns impossible the creation of exact conic curves and surfaces.

This limitation can be circumvented in the case of the resolution of elliptic PDEs using immersed IGA.
As discussed in~\cite{Antolin_2019b}, in those cases it is possible to approximate the geometry of the cut elements $Q\cap\Omega\ \forall Q\in\T^{\Gamma}_h(\Omega)$ by means of B\'ezier curves and surfaces of degree $p$, the same as the solution's discretization, and still preserve optimal approximation properties.
\end{myremark}

\section{Quadrature-free integration of polynomials over B-Reps} \label{sec:integration}
In this section, we deal with the integration of polynomials over a domain~$V$ whose bounding faces~$\Fi$ are represented as trimmed B\'ezier surfaces, as described in the previous section.
More specifically, we seek to compute the integral:
\begin{align} \label{eq:goalIntegral}
  I^{3\textup{D}} = \iiint_{V} a \,\diff V\,,
\end{align}
where $a:V\to\mathbb{R}$ is a polynomial function.
This addresses the computation of the integrals $\mathsf{I}^{3\textup{D}}_{Q,\alpha}$ over cut elements $Q\cap\Omega$ as described in~\eqref{eq:integrals}.

The approach presented in this section consists in the successive application of the divergence theorem, as similarly done, for instance, in~\cite{Lasserre_1998,Duczek_2015,Sudhakar_2014,Antonietti_2018}.
Let us first recall here the classical divergence theorem, also known as Gauss-Ostrogradsky's theorem.
\begin{theorem}
  Let $V$ be a subset of $\mathbb{R}^3$ which is compact and has a piecewise smooth boundary $\partial V$.
  Let $\bm{A}$ be a three-dimensional vector field, such that $\bm{A}:V\to\mathbb{R}^3$ and $\bm{A}\in[C^1(V)]^3$, then:
  \begin{align}\label{eq:thdivergence}
    \iiint_{V} \nabla\cdot \bm{A}\, \diff V = \iint_{\partial V} \bm{A}\cdot\bn\, \diff S\,,
  \end{align}
  where $\nabla\cdot$ is the divergence operator and $\bn:\partial V\to\mathbb{R}^3$ is the outward pointing unit normal on the boundary $\partial V$.
  \label{th:divergence}
\end{theorem}

By applying the divergence theorem, the three-dimensional integral~\eqref{eq:goalIntegral} is transformed into, first, surface, and then line integrals that can be evaluated analytically with machine precision accuracy.
This is possible in the present context due to the polynomial nature of the successive integrands which ease the formation of the antiderivatives involved in the integration process.

\subsection{From volume integral to surface integrals}

In order to apply the divergence theorem, let us first rewrite the initial integral~\eqref{eq:goalIntegral} in the same form as the one in~\eqref{eq:thdivergence}:
\begin{align}
  I^{3\textup{D}} = \iiint_{V} \nabla\cdot\bA \,\diff V\,.
  \label{eq:goalIntegralReform}
\end{align}
The vector field~$\bA:V\to\mathbb{R}^3$ can be expressed as:
\begin{equation}
  \bA(\bx)
  = A_1(\bx)\Vi{e}{1}
  + A_2(\bx)\Vi{e}{2}
  + A_3(\bx)\Vi{e}{3},
\end{equation}
with $\Vi{e}{i}$ as the Cartesian unit vectors and $Q_i:V\to\mathbb{R}$ as the antiderivatives of $a$, computed by:
\begin{equation}
  \label{eq:antiderivative_0}
  \begin{split}
    A_1(x_1,x_2,x_3) &= \alpha_1\int_{0}^{x_1} a(\sigma,x_2,x_3)\diff{\sigma} + \beta_1,\\
    A_2(x_1,x_2,x_3) &= \alpha_2\int_{0}^{x_2} a(x_1,\sigma,x_3)\diff{\sigma} + \beta_2,\\
    A_3(x_1,x_2,x_3) &= \alpha_3\int_{0}^{x_3} a(x_1,x_2,\sigma)\diff{\sigma} + \beta_3.
  \end{split}
\end{equation}
Here $\alpha_1$, $\alpha_2$, $\alpha_3$, $\beta_1$, $\beta_2$, and $\beta_3$ are real constants, such that $\alpha_1+\alpha_2+\alpha_3=1$.
Since $a$ is a polynomial function, the computation of the antiderivatives in~\eqref{eq:antiderivative_0} is straightforward (see Appendix~\ref{sec:appA}).
Furthermore, due to this polynomial nature, the continuity requirements of the divergence theorem are granted for the vector field~$\bA$.

Applying the divergence theorem to~\eqref{eq:goalIntegralReform} we obtain:
\begin{align}
  \label{eq:surfaceIntegral}
  I^{3\textup{D}} = \iint_{\partial V} \bA \cdot \bn \,\diff S\,,
\end{align}
where we recall that $\bn:\partial V\to\mathbb{R}^3$ is the outward pointing unit normal on the boundary $\partial V$.
Recalling the definition of the boundary $\partial V$ in~\eqref{eq:brepPhysical}, the integral~\eqref{eq:surfaceIntegral} can be split as:
\begin{align}
  \label{eq:int_2D_0}
  I^{3\textup{D}}
  = \sumi \,I_i^{2\textup{D}}
  = \sumi \iint_{\Fi} {\bA \cdot \bn_i} \,\diff S_i\,,
\end{align}
where $\bn_i$ are the outward pointing unit normals of the surfaces $\Si$, $i=1,\dots,n_F$.
Exploiting the parametric representation of the surfaces $\Si$, these unit normal vector fields can be expressed as:
\begin{equation}
  \label{eq:normal}
  \bn_i :
  \text{Img}(\Si) \to \mathbb{R}^3\,,~
  \bx
  \mapsto
  \bigg(
  \frac{\bN_i}{\Vert\bN_i\Vert} \circ \Si^{-1}
  \bigg)(\bx)\,,
\end{equation}
where the normal vectors~$\bN_i$ are computed as: 
\begin{equation} \label{eq:Ni}
  \bN_i :
  \operatorname{Dom}(\Si) \to \mathbb{R}^3\,,~
  \bxh
  \mapsto
  \bigg(
    \frac{\partial \Si}{\partial \hat{x}_1} \times \frac{\partial\Si}{\partial\hat{x}_2}
  \bigg)(\bxh)\,.
\end{equation}
In~\eqref{eq:Ni} we assumed that the surface parameterization is oriented such that the cross-product $\bN_i$ points out of $V$.
Plugging~\eqref{eq:normal} into the expression of the surface integrals~$I^{2\textup{D}}$ in~\eqref{eq:int_2D_0}, they become:
\begin{align}
  \label{eq:int_2D_0bis}
  I_i^{2\textup{D}}
  =
  \iint_{\Fi} \bA\cdot \left(\frac{\bN_i}{\Vert\bN_i\Vert} \circ \Si^{-1}\right)  \,\diff S_i\,,
\end{align}
for~$i = 1,\dots{},n_F$.
And pulling back these integrals to the parametric domain of $\Si$, we obtain:
\begin{equation}
  \label{eq:int_2D_1}
  I_i^{2\textup{D}}
  =
  \iint_{\Fhi} \rhi \,\dbxh\,,
\end{equation}
where the integrands $\rhi$ are defined as:
\begin{equation} \label{eq:rhi}
  \begin{split}
    \rhi :
    \operatorname{Dom}(\Si) &\to \mathbb{R}\,,\\
    \bxh
    &\mapsto
    \rhi(\bxh)=
    \big(\bA\circ\Si\big)(\bxh) \cdot \bN_i(\bxh)\,.
  \end{split}
\end{equation}
Interestingly, the normalization and the inversion involved in the definition of the unit normal vectors~\eqref{eq:normal} vanish after the pull-back, as observed in~\cite{Gonzalez_Ochoa_1998}, for instance.
Furthermore, as the surface~$\Si$ is assumed to be polynomial, then the composition $\bA \circ \Si$ is also a polynomial bivariate, but with a higher degree.
Additionally, the non-normalized normal vector field~$\bN_i$ is also a polynomial since it is computed as the product of polynomial terms (the partial derivatives of $\Si$ are polynomials).
Finally, the scalar product of two polynomial vector fields, $\bA\circ\Si$ and $\bN_i$, is a polynomial scalar field.
Consequently, $\rhi$ is a polynomial.
We refer the interested reader to Appendix~\ref{sec:appA} for all the details.

\begin{myremark}
The integrals $I_i^{2\textup{D}}$ in~\eqref{eq:int_2D_1} are equivalent to the boundary integrals $\mathsf{I}^{2\textup{D}}_{Q,\beta}$ depicted in~\eqref{eq:integrals} and required for the assembly of boundary conditions in immersed methods (see Section~\ref{sec:immersed_iga}).
\end{myremark}
\begin{myremark} \label{remark:nontrimmed2D}
  In the case of non-trimmed B\'ezier surfaces, like the one depicted in Figure~\ref{fig:splitPiecewisePolySurf},
  the integrals~\eqref{eq:int_2D_1} can be easily evaluated analytically using Equation~\eqref{eq:multipol:intg}.
\end{myremark}
\begin{myremark} \label{remark:alphas}
  In some situations the normal fields $\bn_i$ of the surfaces $\Si$ may be aligned with one of three the Cartesian axes.
  This occurs quite often in the case of immersed methods for solving PDEs, presented in Section~\ref{sec:immersed_iga}, in which the integration domains $V$ correspond to the cut elements $Q\cap\Omega\ \forall Q\in\T_h(\Omega)$ of the grid embedded in a B-Rep geometry.
  In that particular situation many faces $\Fi$ will be planar trimmed surfaces parallel to the Cartesian axes.
  For those cases, a wise choice of the coefficients $\alpha_1$, $\alpha_2$, and $\alpha_3$ in the antiderivatives~\eqref{eq:antiderivative_0} will make the scalar product $\bA\cdot\bn_i$ vanish, minimizing the number of two-dimensional integrals to be computed.
  For instance, in the case of a face $\Fi$ that is perpendicular to the $z$ Cartesian axis, choosing $\alpha_3=0$ will make the term $\bA\cdot\bn_i$ vanish.
  Nevertheless, for a given domain $V$ the coefficients $\alpha_1$, $\alpha_2$, and $\alpha_3$ must be set once and for all, and cannot be independently chosen for every face $\Fi$ of $V$.
  Thus, an optimal strategy may be to set $\alpha_1$, $\alpha_2$, and $\alpha_3$ independently for every $V$ such that the largest number of surface integrals vanish for that specific domain.
\end{myremark}

\subsection{Evaluating the surface boundary integrals}
Applying again the divergence theorem~\eqref{eq:thdivergence}, we can transform the two-dimensional integrals $I_i^{2\textup{D}}$ in~\eqref{eq:int_2D_1} into line integrals as:
\begin{equation}
  I_i^{2\textup{D}}
  = \int_{\partial \Fhi} \bRhi\cdot\hat{\beme}_i \,\diff \ell_i \,,
  \label{eq:int_1D_0}
\end{equation}
where $\hat{\beme}_i:\partial\Fhi\to\mathbb{R}^2$ is the outward pointing unit normal on the boundary $\partial \Fhi$.
The vector field $\bRhi:\operatorname{Dom}(\Si)\to\mathbb{R}^2$ is defined such that $\rhi = \hnabla\cdot\bRhi$, as for instance:
\begin{equation}
  \begin{split}
    \bRhi(\hat x_1,\hat x_1)
    &= \bigg(\delta_1\int_{0}^{\hat x_1} \rhi(\sigma,\hat x_2) \,\diff\sigma + \epsilon_1 \bigg)\Vi{e}{1}\\
    &+ \bigg(\delta_2\int_{0}^{\hat x_2} \rhi(\hat x_1,\sigma) \,\diff\sigma + \epsilon_2 \bigg)\Vi{e}{2}\,,
  \end{split}
  \label{eq:antiderivative_1}
\end{equation}
and $\delta_1$, $\delta_2$, $\epsilon_1$, and $\epsilon_2$ are real constants, such that $\delta_1+\delta_2=1$.

Splitting the boundary $\partial\Fhi$ according to~\eqref{eq:Fhi} we obtain:
\begin{equation}
  I_i^{2\textup{D}}
  =
  \sumj I_{i,j}^{1\textup{D}}
  =
  \sumj \int_{\ghij} \bRhi\cdot\hat{\beme}_{i,j} \,\diff \ell_{i,j}\,,
  \label{eq:int_1D_1}
\end{equation}
where $\hat{\beme}_{i,j}:\text{Img}(\chij)\to\mathbb{R}^2$ are the outward pointing unit normals of the curves $\chij$, $i=1,\dots,n_{c,i}$.
Exploiting the parametric representation of the curves $\chij$, these unit normal vector fields can be expressed as:
\begin{equation}
  \label{eq:normal2D}
  \hat{\beme}_{i,j} :
  \text{Img}(\chij) \to \mathbb{R}^2\,,~
  \bxh
  \mapsto
  \bigg(
  \frac{\hat{\bM}_{i,j}}{\Vert\hat{\bM}_{i,j}\Vert} \circ \chij^{-1}
  \bigg)(\bxh)\,.
\end{equation}
where the normal vectors~$\hat{\bM}_{i,j}$ are computed as: 
\begin{equation} \label{eq:Mi}
  \hat{\bM}_{i,j} :
  \operatorname{Dom}(\chij) \to \mathbb{R}^2\,,~
  \xt
  \mapsto
  \frac{\diff\chij}{\diff\xt}(\xt) \times \Vi{e}{3}
  \,.
\end{equation}
In the previous expression we assume that the curves $\chij$ are oriented such as the external boundaries of $\Fhi$ present a counter-clockwise orientation, while the internal ones are clockwise oriented (see Figure~\ref{fig:surfaceTOcurves}).

Plugging~\eqref{eq:normal2D} into the expression of the line integrals~$I^{1\textup{D}}$ involved in~\eqref{eq:int_1D_1}, they become:
\begin{align}
  \label{eq:int_1D_2}
  I_{i,j}^{1\textup{D}}
  =
  \int_{\ghij} \bRhi\cdot \bigg(\frac{\hat{\bM}_{i,j}}{\Vert\hat{\bM}_{i,j}\Vert} \circ \chij^{-1}\bigg)  \,\diff \ell_{i,j}\,.
\end{align}
Finally, pulling back these integrals to the parametric domain of the underlying curves $\chij$, we obtain:
\begin{equation}
  \label{eq:int_1D_3}
  I_{i,j}^{1\textup{D}}
  =
  \int_{0}^{1} \big(\bRhi\circ\chij\big) \cdot \hat{\bM}_{i,j} \,\dxt\,,
\end{equation}
where, as for the two-dimensional case, the normalization and the inversion involved in the definition of the unit normal vectors~\eqref{eq:normal2D} vanish after the pull-back.
We gather all the integrand terms together as:
\begin{equation}
  \label{eq:int_1D_4}
  I_{i,j}^{1\textup{D}}
  =
  \int_{0}^{1} \tilde{t}_{i,j} \,\dxt\,,
\end{equation}
where
\begin{equation} \label{eq:tij}
  \begin{split}
    \tilde{t}_{i,j} :
    &\operatorname{Dom}(\chij)=[0,1] \to \mathbb{R}\,,\\
    &\xt
    \mapsto
    \tilde{t}(\xt)=
    \big(\bRhi\circ\chij\big)(\xt) \cdot \hat{\bM}_{i,j}(\xt)\,.
  \end{split}
\end{equation}
As the curve~$\chij$ is a B\'ezier, the composition $\bRhi\circ\chij$ is a higher degree univariate polynomial.
Additionally, the non-normalized normal vector field~$\hat{\bM}_{i,j}$ is also a polynomial since it is computed from B\'ezier derivatives.
Finally, the scalar product of two polynomial vector fields, $\bRhi\circ\chij$ and $\hat{\bM}_{i,j}$, is a polynomial scalar field.
Consequently, $\tilde{t}_{i,j}$ is a polynomial.
Therefore, the integrals~\eqref{eq:int_1D_4} can be easily evaluated in an analytic way (see Equation~\eqref{eq:pol:intg}),  with machine precision accuracy, without the need of quadrature schemes.
Further details for the case of Bernstein polynomials are provided in Appendix~\ref{sec:appA}.

\begin{myremark} \label{remark:deltas}
  The Remark \ref{remark:alphas} is extensible to the line integrals detailed above.
  In some situations (see for instance Figure~\ref{fig:surfaceTOcurves}), some boundaries $\ghij$ may be aligned with the Cartesian axes.
  In those cases, the constants $\delta_1$ and $\delta_2$ arising in the antiderivatives~\eqref{eq:antiderivative_1} can be chosen such as the product $\bRhi\cdot\hat{\beme}_{i,j}$ vanishes in some of those boundaries.
  These constants can be chosen independently for every face integral $I_{i}^{2\textup{D}}$ such as the number of 1D integrals to be evaluated is minimized.
\end{myremark}

\begin{figure*}[tp]
  \centering
  \includegraphics[width=\linewidth]{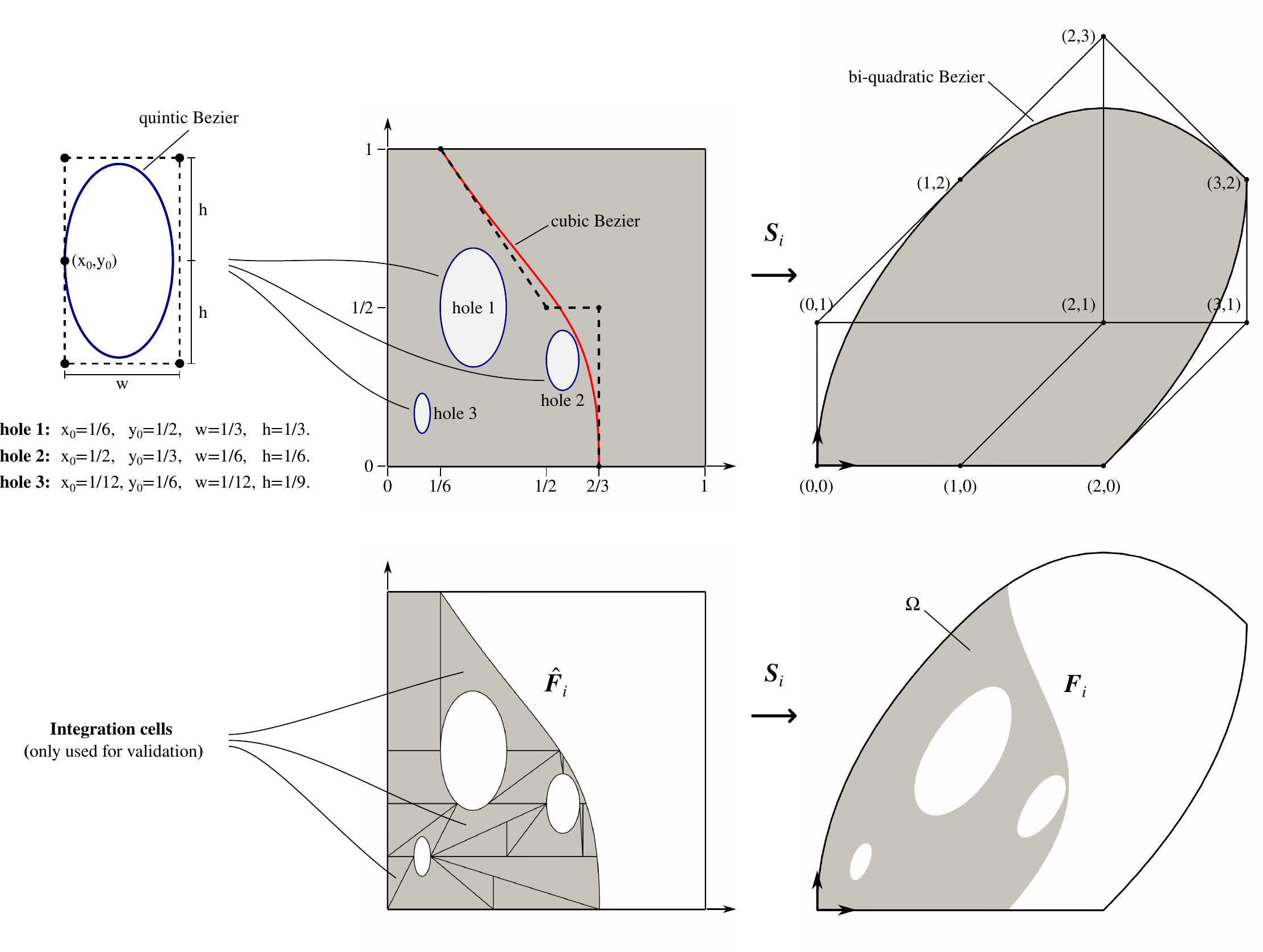}
  \caption{The two-dimensional trimmed geometry for the validation of the quadrature-free integration procedure.}
  \label{fig:integration1}
\end{figure*}

\subsection{Polynomial degree} \label{sec:degrees}
The reader may have notice that due to the involved compositions, $\bA\circ\Si$ and $\bRhi\circ\chij$, as well as the products of B\'ezier curves and surfaces, the resulting polynomial term $\tilde{t}_{i,j}$ can potentially present a very high degree.
In this section we detail the computation of this degree, as well as the order of other terms involved in the intermediate steps.

For the sake of simplicity, hereinafter we assume that the polynomial $a$ to integrate, as well as the B\'ezier mappings $\Si$ and $\chij$, have constant degrees along all their parametric directions and for all their components:
\begin{equation} \label{eq:degrees0}
  a     \in\mathbb{Q}_{r,r,r};~
  \Si   \in\mathbb{Q}_{s,s}\times\mathbb{Q}_{s,s}\times\mathbb{Q}_{s,s};~
  \chij \in\mathbb{Q}_c\times\mathbb{Q}_c,
\end{equation}
with $r\geq0$, $s>0$, and $c>0$, and where the polynomial spaces $\mathbb{Q}$ follow the notation introduced in Section~\ref{sec:polapprox}.
According to the definitions~\eqref{eq:Ni} and \eqref{eq:Mi} it is straightforward to obtain the degrees of the fields $\bN_i$ and $\hat\bM_{i,j}$ as:
\begin{equation}
  \begin{split}
    \bN_i &\in\mathbb{Q}_{2s-1,\,2s-1}\times\mathbb{Q}_{2s-1,\,2s-1}\times\mathbb{Q}_{2s-1,\,2s-1}\,,\\
    \hat\bM_{i,j} &\in\mathbb{Q}_{c-1}\times\mathbb{Q}_{c-1}\,,
  \end{split}
\end{equation}
and using~\eqref{eq:antiderivative_0}, the order of $\bA$ is computed as:
\begin{equation}
  \bA \in\mathbb{Q}_{r+1,\,r,\,r}\times\mathbb{Q}_{r,\,r+1,\,r}\times\mathbb{Q}_{r,\,r,\,r+1}\,.
\end{equation}
Thus, the degrees of $\bA\circ\Si$ and $\rhi$ (recall Equation~\eqref{eq:rhi}) are:
\begin{equation}
  \begin{split}
    &\bA\circ\Si \in\mathbb{Q}_{t,\,t}\times\mathbb{Q}_{t,\,t}\times\mathbb{Q}_{t,\,t},\quad{}t= 2\left(3 r + 1\right),\\
    &\rhi \in\mathbb{Q}_{3s\left(r+1\right)-1,\,3s\left(r+1\right)-1}\,.
\end{split}
\end{equation}

\begin{figure*}[t]
  \centering
  \includegraphics[width=0.9\linewidth]{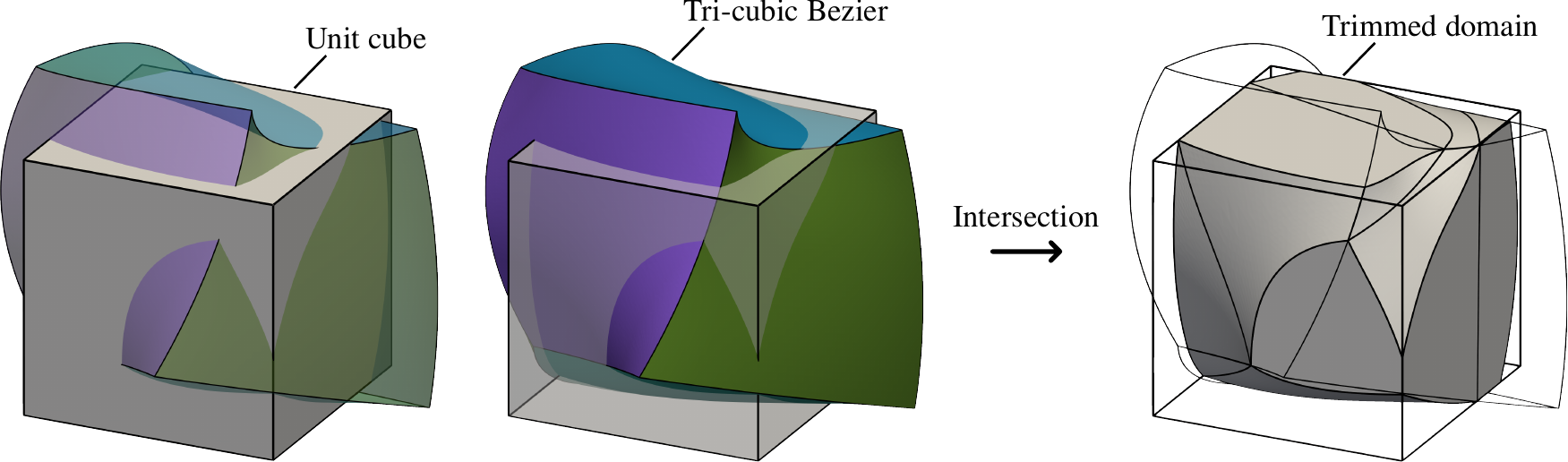}
  \caption{The three-dimensional trimmed geometry for the validation of the quadrature-free integration procedure.}
  \label{fig:integration2}
\end{figure*}

Analogously to the case of $\bA$, the degree of $\bRhi$ (Equation~\ref{eq:antiderivative_1}), and its composition $\bRhi\circ\chij$, are simply computed as:
\begin{equation}
  \begin{split}
    &\bRhi \in\mathbb{Q}_{3s\left(r+1\right),\,3s\left(r+1\right)-1}\times\mathbb{Q}_{3s\left(r+1\right)-1,\,3s\left(r+1\right)}\,,\\
    &\bRhi\circ\chij \in\mathbb{Q}_{6sc\left(r+1\right)-c}\times\mathbb{Q}_{6sc\left(r+1\right)-c}\,.
  \end{split}
\end{equation}
Finally, the polynomial term $\tilde{t}_{i,j}$ presents a degree:
\begin{equation} \label{eq:final_degree}
  \tilde{t}_{i,j} \in\mathbb{Q}_{6sc\left(r+1\right)-1}\,.
\end{equation}

The degree of $\tilde{t}_{i,j}$ can be potentially very high what may induce numerical instabilities.
Nevertheless, in the examples of Section~\ref{sec:poisson3Dsimple} very high order polynomials were involved (in the order of hundreds) but no instabilities were noticed.
This is due to the fact that we use B\'ezier curves and surfaces that are expressed in terms of Bernstein polynomials, known to be numerically more stable than other choices, as, for instance, monomial or Lagrange bases.
Along this work, we compute derivatives, integrals, additions, and multiplications of Bernstein polynomials, that are stable operations, but we never evaluate polynomials.
See Appendix~\ref{sec:appA} for further details.

\section{Numerical experiments}
\label{sec:numericalExamples}

In this section we show the performance of the presented quadrature-free approach by means of numerical experiments.
In a first set of examples, in Section~\ref{sec:num_integrals}, we apply the method to the computation of simple integrals in 2D and 3D domains and compare them with standard methods based on the use of boundary-conforming quadrature schemes.
Afterwards, in Section~\ref{sec:num_poisson} we apply it to the solution of elliptic PDEs using the immersed isogeometric framework presented in Section~\ref{sec:immersed_iga}.

\subsection{Computation of integrals over B-reps} \label{sec:num_integrals}

\newcommand{\mass}{\ensuremath{M}}
\newcommand{\centermass}{\ensuremath{\bm{C}_M}}

Figures~\ref{fig:integration1} and~\ref{fig:integration2} present two numerical studies used to validate the presented integration strategy.
The two-dimensional case, described in Figure~\ref{fig:integration1}, consists in a quadratic B\'ezier surface which is trimmed by three holes and a vertical curved slice.
The three-dimensional case, described in Figure~\ref{fig:integration2}, involves a trimmed domain defined by the intersection of a cube and a free-form cubic trivariate.
We compute the mass~\mass{} and the center of gravity~\centermass{} of these two geometries, defined by:
\begin{subequations} \label{eq:mass}
  \begin{align}
  \mass &= \int_{V} \rho(\bx) \dbx,\\
  \centermass &= \frac{1}{\mass} \int_{V} \bx\rho(\bx) \dbx\,,
  \end{align}
\end{subequations}
where the density is considered to be constant $\rho=1$.

\begin{table*}[t]
  \centering
  \begin{tabular}{@{}lrrr@{}}
    \toprule
    & \multicolumn{1}{c}{Reference} & \multicolumn{1}{c}{Quad-free} & \multicolumn{1}{c}{Relative diff.} \\ \midrule
    2D geo: \mass{}                             & 2.100230243261870         & 2.100230243261870         & $<10^{-15}$              \\
    \phantom{2D geo: }$\centermass{}\cdot\bm{e}_1$ & 0.914136125211735         & 0.914136125211735         & $<10^{-15}$              \\
    \phantom{2D geo: }$\centermass{}\cdot\bm{e}_2$ & 0.859802811586580         & 0.859802811586580         & $<10^{-15}$
    \\[5pt]
    3D geo: \mass{}                             & 0.444790448933688         & 0.444790378608127         & $1.58\times{}10^{-7}$    \\
    \phantom{3D geo: }$\centermass{}\cdot\bm{e}_1$ & 0.469169723257000         & 0.469169674580198         & $1.03\times{}10^{-7}$    \\
    \phantom{3D geo: }$\centermass{}\cdot\bm{e}_2$ & 0.400642146493445         & 0.400642138814180         & $1.91\times{}10^{-8}$    \\
    \phantom{3D geo: }$\centermass{}\cdot\bm{e}_3$ & 0.457115007608867         & 0.457114990479802         & $3.74\times{}10^{-8}$    \\   \bottomrule
  \end{tabular}
  \caption{Comparison of the quadrature-free integration for the 2D and 3D trimmed geometries depicted in Figures~\ref{fig:integration1} and~\ref{fig:integration2}, respectively.
           The mass and the center of mass are evaluated and compared to reference values obtained with an alternative approach based on reparameterization.}
  \label{tab:integration}
\end{table*}
Reference values of~\eqref{eq:mass} are obtained through boundary-conformal quadrature schemes created by reparameterizing the interior of $V$ with a technique similar to the one presented in~\cite{Massarwi_2019}.
This approach subdivides the domain of integration and leads to integration sub-cells.
Standard quadrature rules can then be used to integrate numerically.
For the sake of comparison, an overkill number of quadrature points were used within each integration cell for both examples.

The obtained results are presented in Table~\ref{tab:integration}.
For the 2D-geometry (Figure~\ref{fig:integration1}), the computed relative differences, compared with the reparameterization approach, are below $10^{-15}$, \emph{i.e.}, close to machine precision.
Nevertheless, for the 3D-geometry (Figure~\ref{fig:integration2}), relative differences of the order of $10^{-7}$ were noticed.

\begin{myremark} \label{remark:tolerance}
We associate the larger differences in the 3D case to the intrinsic tolerances involved in some geometric operations.
In this work we employ algorithms provided by Open CASCADE Technology~\cite{OpenCASCADE} which is an open source \texttt{C++} library designed for geometric modeling applications.
For instance, in the specific case of surface-surface intersections between B-spline or B\'ezier surfaces, Open CASCADE limits the lowest tolerance to $10^{-7}$, what truncates the achievable accuracy and agrees with the results reported in Table~\ref{tab:integration}.
Similar tolerances apply to other non-linear operations.
These limitations are not exclusive of Open CASCADE, as similar issues can be found in other commercial and non-commercial geometric kernels available:
Tolerances of the order of $10^{-7}$ are more than enough for most of the applications these tools are designed for.
On the other hand, we use Irit~\cite{Irit}, an open source geometric modeler, for other 2D operations, as it is the case of the computation of intersections between planar spline curves.
The involved tolerances in Irit can be tuned according to our needs, what allows us to reach a higher accuracy for the 2D problem.
In addition, it is important to remark that these limitations pollute the geometrical approximation not just for the presented quadrature-free method, but as well for other approaches, as for instance, for surface and volumetric untrimming, as previously discussed in~\cite{Antolin_2019b}.
Nevertheless, we believe that the obtained results confirm the viability of the quadrature-free integration strategy for 3D~geometries.
\end{myremark}

\begin{myremark}
  For computing the quantities~\eqref{eq:mass} in the case of the 2D-geometry (Figure~\ref{fig:integration1}), the integration procedure can be directly started from Equation~\eqref{eq:int_2D_1}, by replacing $\hat{r}_i(\bxh)$ with
  $\big(\rho\circ\Si\big)(\bxh)$ and $\big(\rho\circ\Si\big)(\bxh)\,\Si(\bxh)\cdot\bm{e}_k,~k=1,2,3$, respectively.
\end{myremark}

\subsection{Immersed isogeometric analysis} \label{sec:num_poisson}

In this section we demonstrate the effectiveness of the quadrature-free approach for solving PDEs in the context of the immersed isogeometric framework presented in Section~\ref{sec:immersed_iga}.
In particular, we perform a series convergence analyses for the Poisson's problem in different 2D (Section~\ref{sec:poisson2D}) and 3D (Section~\ref{sec:poisson3Dsimple}) immersed domains.
Optimal error convergence rates are retrieved in all the cases.
Finally, in Section~\ref{sec:poisson3Dcomplex}, the flexibility and robustness of the proposed approach is demonstrated in the case of geometries that present a level complexity analogous to the ones found in real industrial applications.

For all the studied cases, we consider the approximated Poisson's problem~\eqref{eq:weakpbapprox}, previously discussed in Section~\ref{sec:immersed_iga}.
We adopt manufactured solutions:
\begin{equation} \label{eq:poisson_exact_solutions}
  \begin{split}
    u_{\text{ex}}(x,y) &= \sin(\pi x)\sin(\pi y)&~~\text{in 2D}\,,\\
    u_{\text{ex}}(x,y,z) &= \sin(\pi x)\sin(\pi y)\sin(\pi z)&~~\text{in 3D}\,,
  \end{split}
\end{equation}
except for the complex geometries in Section~\ref{sec:poisson3Dcomplex}.
Accordingly, the source and Neumann terms, $f$ and $g$, are defined as:
\begin{subequations}
  \begin{align}
    f &= -\Delta{u_{\text{ex}}}\,,\\
    g &= \nabla{u_{\text{ex}}}\cdot\bn\,.
  \end{align}
\end{subequations}
The Dirichlet boundary $\Gamma_D$ will be defined for each particular case, and, consequently, Neumann boundary conditions will be applied on $\Gamma_N=\partial\Omega\setminus\Gamma_D$.

The choice of such regular functions as target solutions (Equation~\eqref{eq:poisson_exact_solutions}) is motivated by the aim of focusing our study on the consistency error, mainly controlled by numerical integration and geometric representation errors, while keeping the discretization error small.
The approximation properties of trimmed spline spaces for the solution of elliptic PDEs have been previously studied in \cite{Antolin_2019b}.

\subsubsection{Poisson's problem for 2D trimmed-geometries} \label{sec:poisson2D}

\begin{figure*}[htp]
  \centering
  \includegraphics[width=0.9\linewidth]{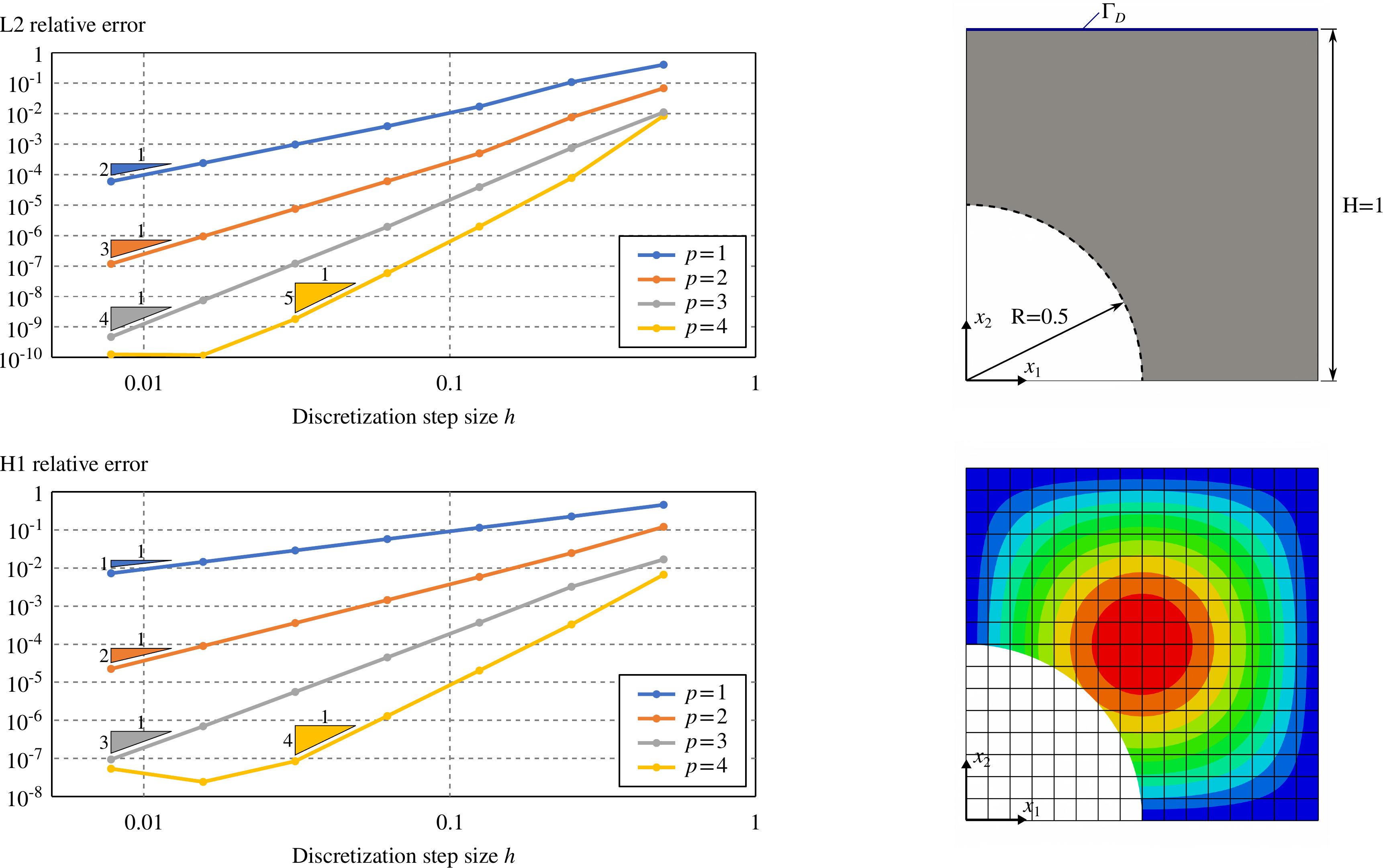}
  \caption{Poisson's problem over a square with a circular trimmed region.}
  \label{fig:poisson2Da}
\end{figure*}

\begin{figure*}[htp]
  \centering
  \includegraphics[width=0.9\linewidth]{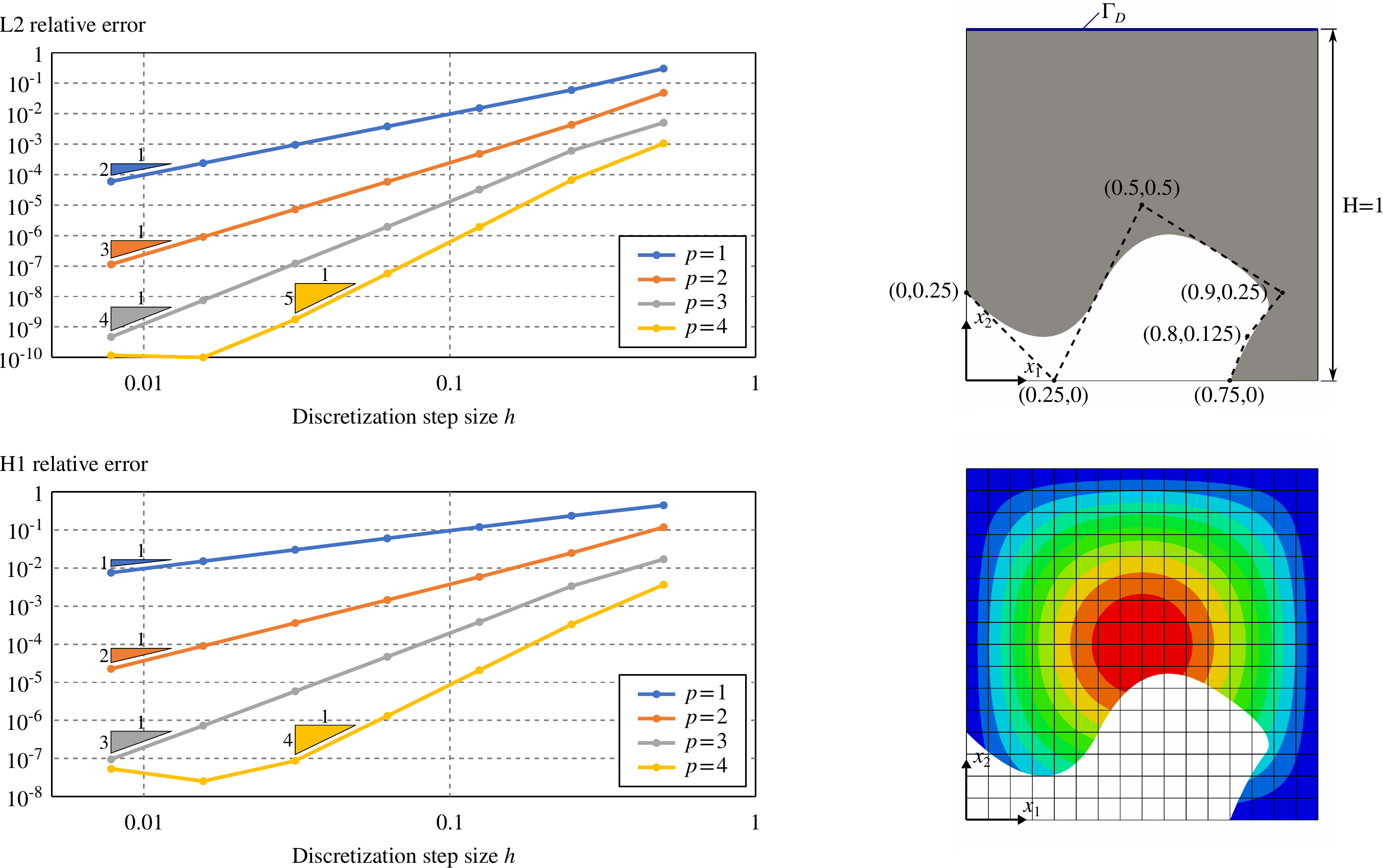}
  \caption{Poisson's problem over a square with a free-form trimmed region.}
  \label{fig:poisson2Db}
\end{figure*}

\begin{figure*}[htp]
  \centering
  \includegraphics[width=\linewidth]{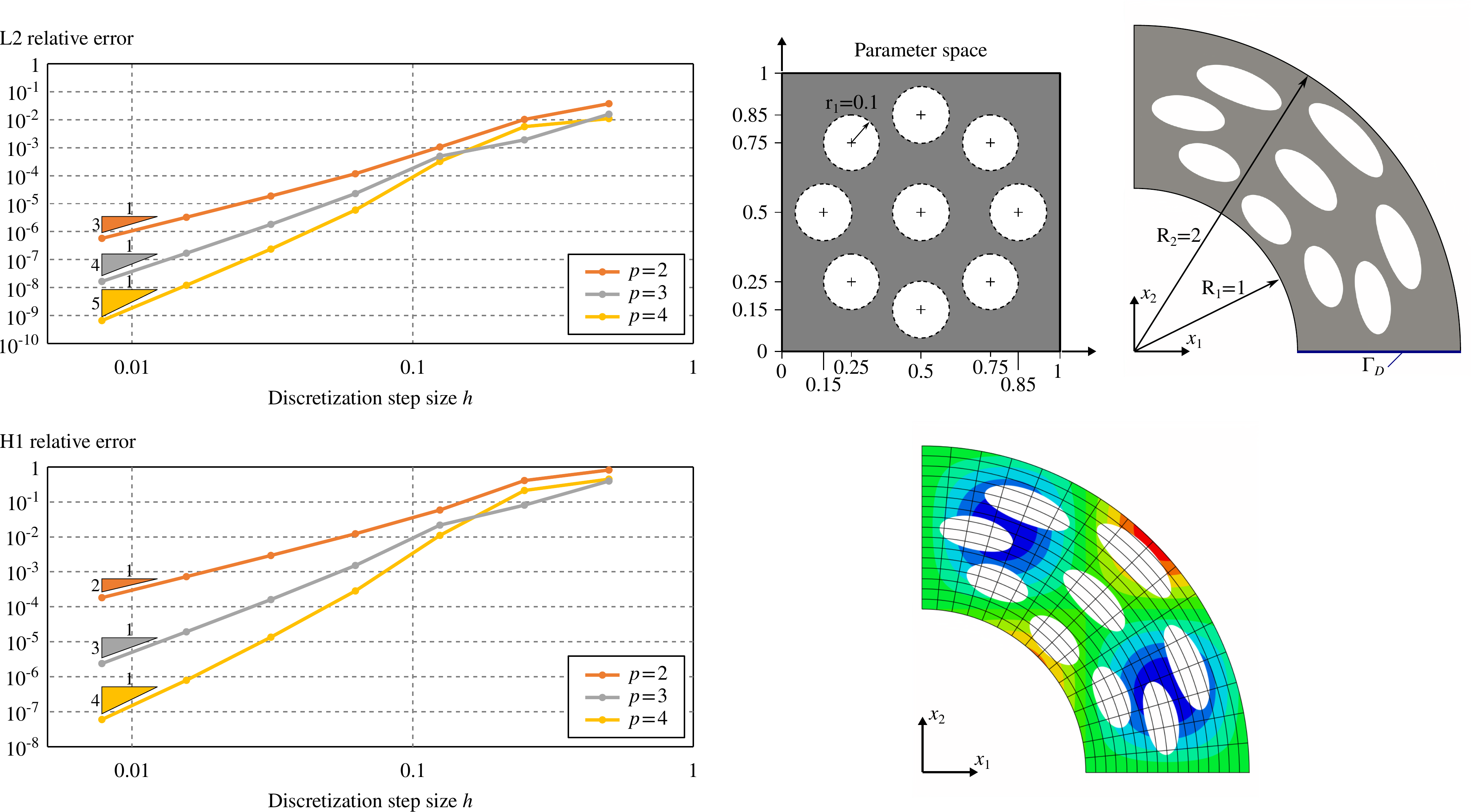}
  \caption{Poisson's problem over a one-quarter annulus with several holes.}
  \label{fig:poisson2Dc}
\end{figure*}

Let us firstly tackle the Poisson's problem for several two-dimensional problems:
\begin{itemize}
\item a square with a circular hole (Figure~\ref{fig:poisson2Da}),
\item a square with a free-form hole (Figure~\ref{fig:poisson2Db}),
\item a multi-perforated quarter annulus (Figure~\ref{fig:poisson2Dc}).
\end{itemize}
Several solution degrees are considered: \emph{i.e.}, from $p=1$ for the trimmed squares, and $p=2$ for the annulus, to $p=4$.
Importantly, the presence of conic sections require to perform some geometric approximations such that the integrals in the finite element operators involve only non-rational polynomials.
As already discussed in Remark~\ref{rmk:rationals}, to do so we rely on the results proven in~\cite{Antolin_2019b} which reveal that approximating the elements' geometry using degree~$p$ leads to optimal numerical results.
Therefore, B\'eziers of degree~$p$ are used to approximate the rational geometrical quantities at the element level.

In addition, it is important to remark the presence of a non-identity mapping in the problem depicted in Figure~\ref{fig:poisson2Dc}.
This leads to the introduction of an extra non-polynomial term in the bilinear form (see Remark~\ref{rmk:mapping}) that is approximated through a local polynomial projection, as discussed in Section~\ref{sec:polapprox}.

The solution errors in both $H^1$ and $L^2$ relative norms are evaluated along with the analyses.
The evaluation of these errors is done through the use of the reparameterization approach already employed during the validation of the integrals computed in Section~\ref{sec:num_integrals}.
Optimal convergence rates, $p$ and $p+1$, respectively, are retrieve for the three cases, see again Figures~\ref{fig:poisson2Da}, ~\ref{fig:poisson2Db}, and~\ref{fig:poisson2Dc}.
The numerical solutions obtained with the quadrature-free approach enable to validate the present methodology for two-dimensional cases.

Nevertheless, it is important to remark that for the finest discretizations in the case $p=4$, the error reaches a \emph{plateau} (around $10^{-10}$ for the relative $L^2$ error norms).
For those cases, the discretization error becomes lower than the error induced by geometrical operations as, for instance, the slicing of the domain $\Omega$ into elements.
See the related discussion in Remark~\ref{remark:tolerance}.
Similar \emph{plateaux} were observed in \cite{Antolin_2019b,antolin2021overlapping}.

\subsubsection{Poisson's problem for simple 3D trimmed-geometries} \label{sec:poisson3Dsimple}

\begin{figure*}[tp]
  \centering
  \includegraphics[width=0.8\linewidth]{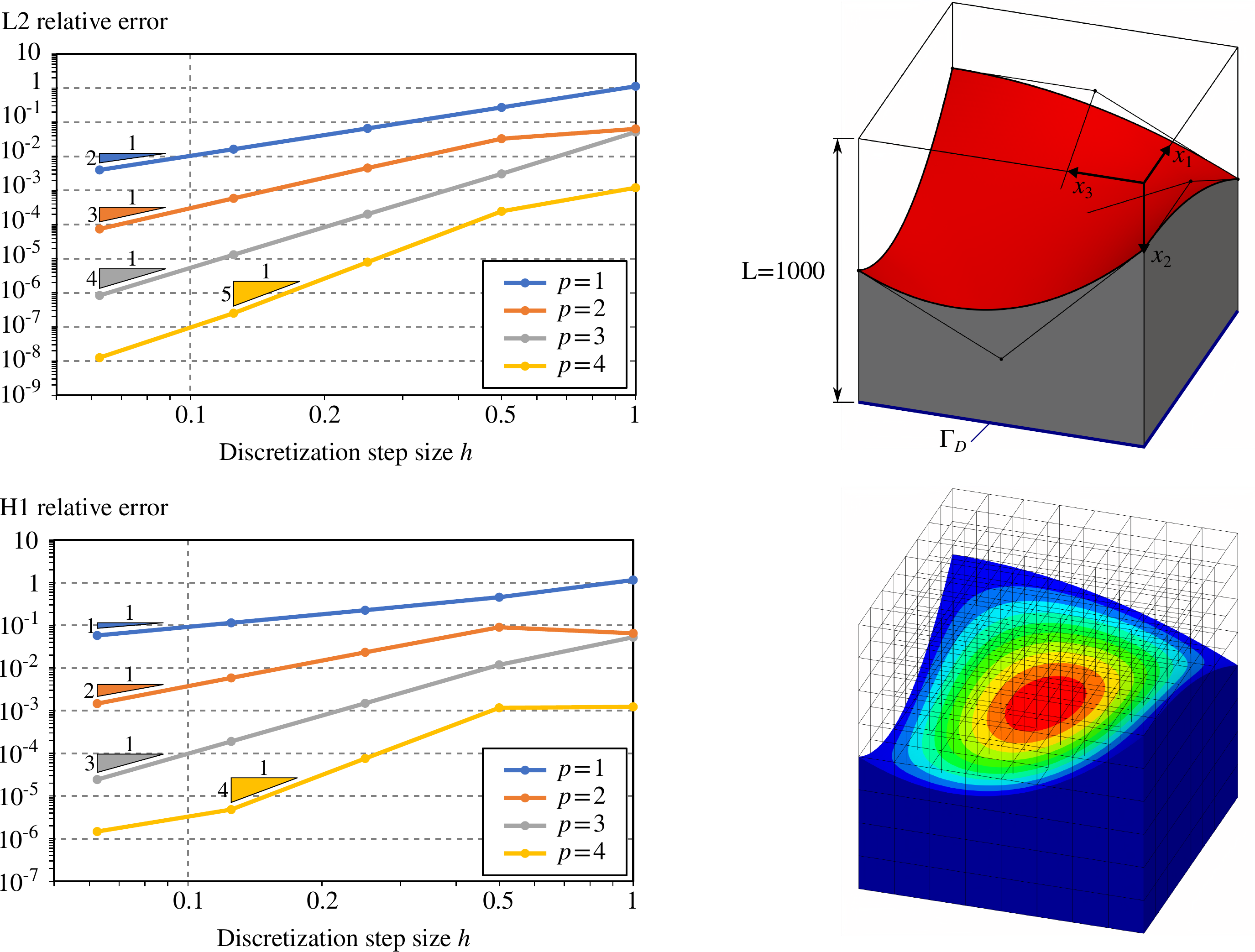}
  \caption{Poisson's problem over a cube with a free-form trimmed region.}
  \label{fig:poisson3Da}
\end{figure*}

\begin{figure*}[tp]
  \centering
  \includegraphics[width=0.8\linewidth]{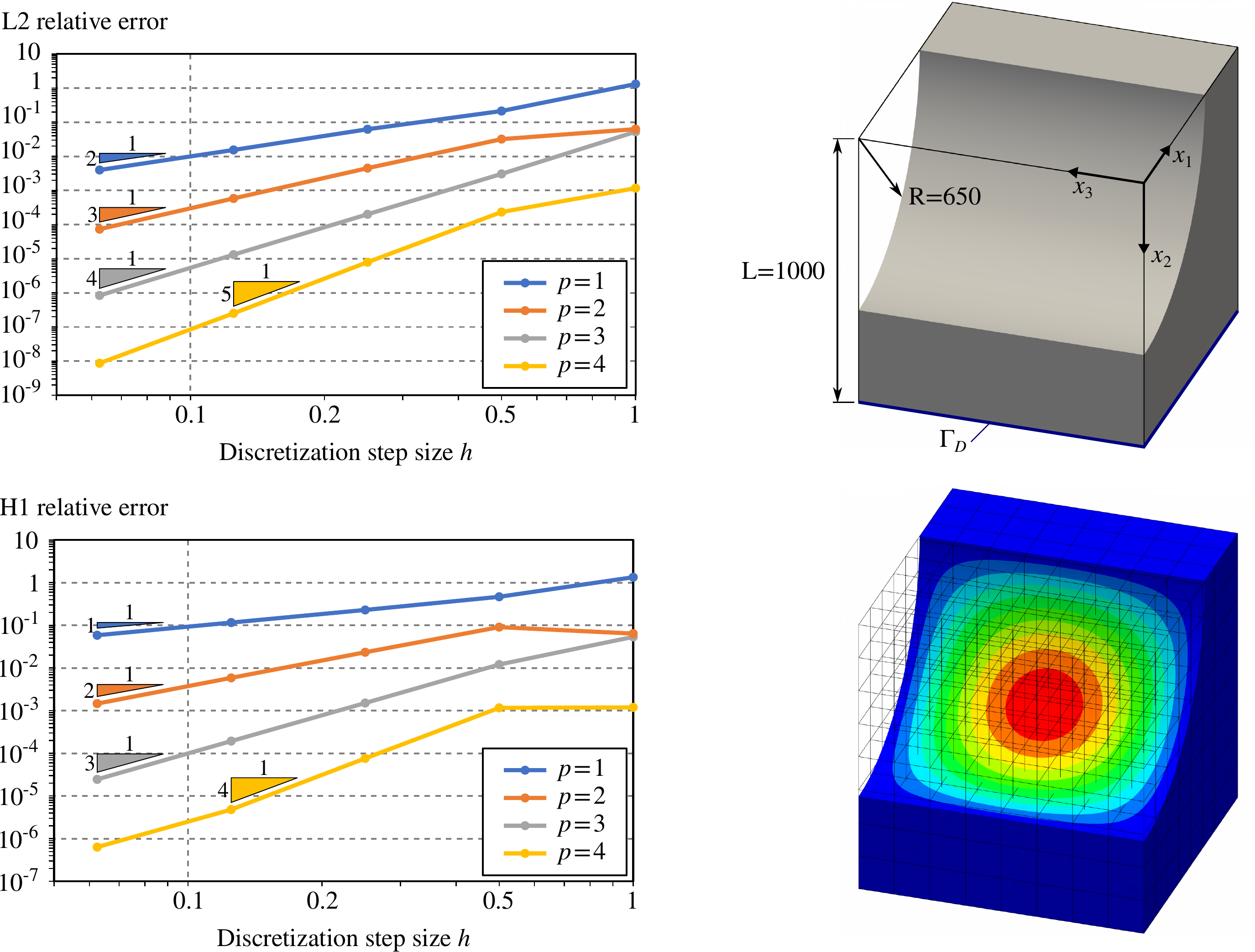}
  \caption{Poisson's problem over a cube with a cylindrical trimmed region.}
  \label{fig:poisson3Db}
\end{figure*}

\begin{figure*}[tp]
  \centering
  \includegraphics[width=0.8\linewidth]{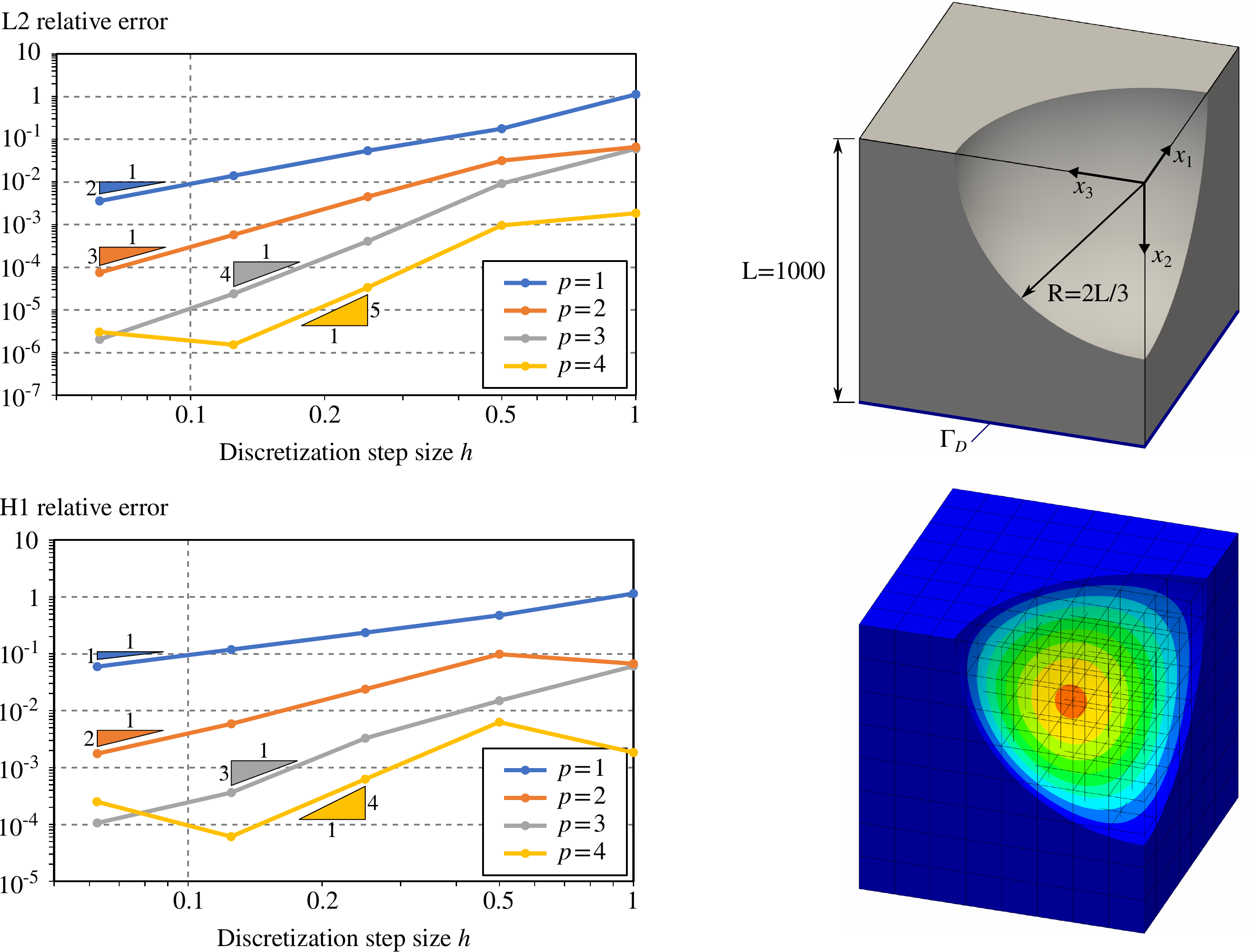}
  \caption{Poisson's problem over a cube with a spherical trimmed region.}
  \label{fig:poisson3Dc}
\end{figure*}

In order to go one step further, we perform several analyses on three-dimensional trimmed domains.
We consider again three trimmed domains.
Each of them consists in a cube with length~$L=1000$ with different trimmed regions:
\begin{itemize}
\item a free-form cut (Figure~\ref{fig:poisson3Da}) which is defined by a bi-quadratic surface with the following control points:
  {\small
    \begin{equation*}
      \mathbf{P}=
      \begin{pmatrix}
        0 & L/2 & L & 0 & L/2 & L & 0 & L/2 & L \\
        L/4 & L/4 & L/2 & 3L/4 & L/2 & L/4 & L/2 & 3L/4 & L/4 \\
        0 & 0 & 0 & L/2 & L/2 & L/2 & L/2 & 3L/4 & L/4 
      \end{pmatrix}
    \end{equation*}%
  }
\item one-quarter of a cylinder (Figure~\ref{fig:poisson3Db}),
\item one-eighth of a sphere (Figure~\ref{fig:poisson3Dc}).
\end{itemize}
As for the 2D-cases, we study the convergence rate in both $H^1$ and $L^2$ relative norms for several spline degrees.
The norms are again evaluated via a reparameterization procedure.
The obtained results confirm the theoretical expectations: Optimal convergence rates are confirmed.

As in the 2D case, for the finest discretization in the case $p=4$ the optimal convergence rate starts to deteriorate.
Again, this is due to the consistency errors introduced by the involved geometric operations, as discussed in Remark~\ref{remark:tolerance}.
Similar results were previously observed in~\cite{Antolin_2019b}.

The curves created by Open CASCADE~\cite{OpenCASCADE} during the surface-surface intersections are represented as B-splines of high degree (for instance, degree 8 for the spherical removal in Figure~\ref{fig:poisson3Dc}) and possibly rational.
Such high order curves may lead to very high degrees during the polynomial compositions, as detailed in Section~\ref{sec:degrees}.
As for the 2D-cases, and according to Remark~\ref{rmk:rationals}, it is always possible to approximate at the element level those geometrical entities with B\'eziers of degree equal to the solution degree.
What turns to be mandatory in the case of rational curves and surfaces.

In all the three numerical examples included in this section the curves arising from surface-surface intersections were approximated at element level using B\'ezier curves with degree $p$.
In the same way, for the cases in Figures~\ref{fig:poisson3Db} and \ref{fig:poisson3Dc}, the underlying rational surfaces were also approximated at element level with B\'ezier surfaces of degree $p$ along both parametric directions.

Let us now study the involved polynomial degrees for the three examples included in this section according to the estimation detailed in Section~\ref{sec:degrees}.
Applying the quadrature-free approach to solve the Poisson's problem~\eqref{eq:weakpbapprox}, we can identify the polynomial integrand $a$ (recall Equation~\eqref{eq:goalIntegral}) with the term $\left.B^{\mathbf{0}}_k \big(\nabla{N_i}\otimes\nabla{N_j}\big)\right\rvert_{Q}\in\mathbb{Q}_{2p,\,2p,\,2p}$ (Equation~\eqref{eq:termsdegrees}, where we assumed $\bar{\diffusivityCoef}$ to be the identity and therefore the projection degrees to be $\mathbf{q}=(0,0,0)$).

Considering, as discussed above, that the degrees of approximated surfaces and curves are $s=c=p$, the final degree of the polynomial term $\tilde{t}_{i,j}$ becomes (recall Equation~\eqref{eq:final_degree}):
\begin{equation}
  \tilde{t}_{i,j} \in\mathbb{Q}_{w}\,,\text{with }w=12 p^3 + 6 p^2 - 1\,.
\end{equation}
Unsurprisingly, the degree $w$ is very high: $w=\left\lbrace 17,\,119,\,377,\,864\right\rbrace$ for $p=\left\lbrace 1,\,2,\,3,\,4\right\rbrace$, respectively.
Nevertheless, despite these high orders, no instabilities were noticed in the results of Figures~\ref{fig:poisson3Da}-\ref{fig:poisson3Dc}.
As previously discussed in Section~\ref{sec:degrees}, this is due to the fact that the proposed integration strategy does not require polynomial evaluations.
An in-depth discussion can be found in Appendix~\ref{sec:appA}.

\subsubsection{Poisson's problem on complex 3D trimmed-geometries} \label{sec:poisson3Dcomplex}

\begin{figure*}[tp]
  \centering
  \includegraphics[width=\linewidth]{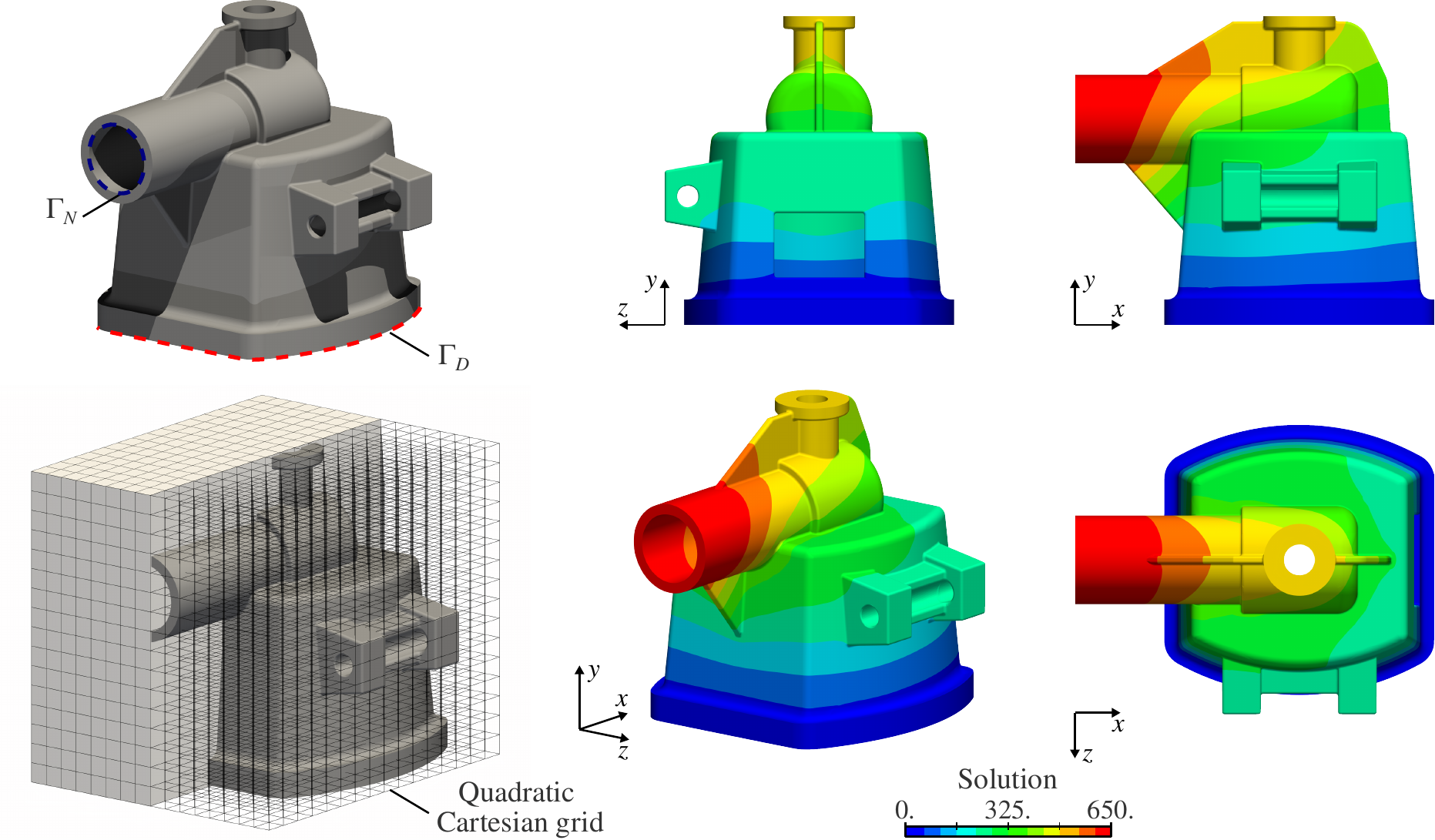}
  \caption{Poisson's problem over the first complex CAD geometry.}
  \label{fig:complex_geom_1}
\end{figure*}

\begin{figure*}[tp]
  \centering
  \includegraphics[width=\linewidth]{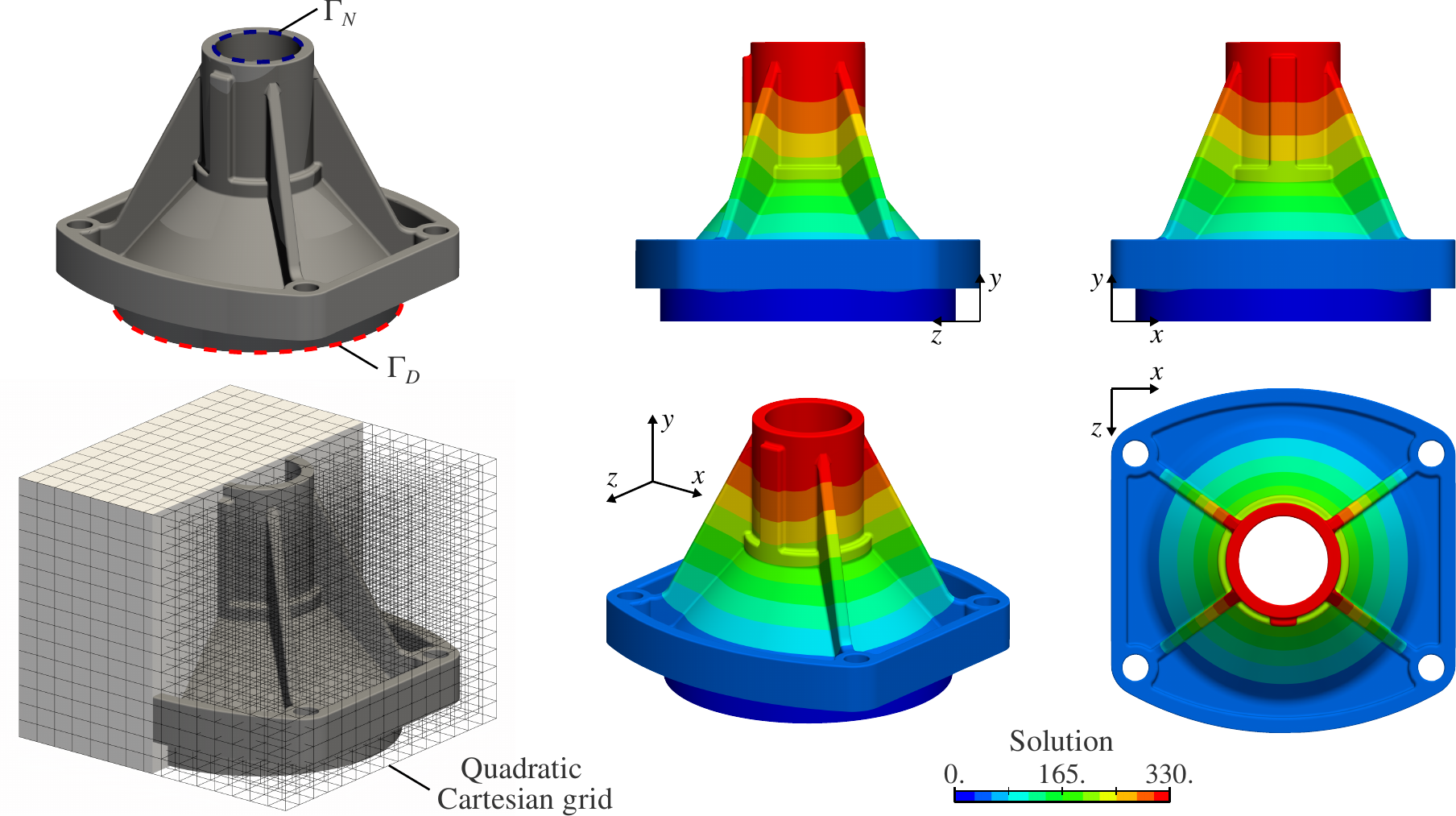}
  \caption{Poisson's problem over the second complex CAD geometry.}
  \label{fig:complex_geom_2}
\end{figure*}

In order to show the viability of the quadrature-free approach to handle complex 3D geometries, we consider the two CAD models shown in Figures~\ref{fig:complex_geom_1} and~\ref{fig:complex_geom_2}.
These B-Rep geometries have been extracted from the Open CASCADE database~\cite{OpenCASCADE}.
Generating a boundary-conforming volumetric parameterization of these geometries is far from a simple task.
Instead, the B-Rep models are immersed into Cartesian grids (see Section~\ref{sec:immersed_iga}).
The solutions are discretized with $C^1$-continuous quadratic B-spline basis functions.
Again, we solve Poisson's problem with homogeneous Dirichlet boundary condition applied on the bottom surfaces and a constant Neumann boundary condition inside the cylindrical tubes (see again Figures~\ref{fig:complex_geom_1} and~\ref{fig:complex_geom_2}).
In order to build the finite element operators, the presented quadrature-free approach is applied.
The obtained solutions are depicted in Figures~\ref{fig:complex_geom_1} and~\ref{fig:complex_geom_2}.
We believe that these two complex geometries highlight the viability of the developed approach to deal with designs of industrial complexity level.

\section{Conclusions} \label{sec:conclusions}

We have presented a novel approach for the solution of partial differential equations on B-Rep geometries by means of immersed isogeometric discretizations that do not require quadrature schemes.
For such purpose, we developed a new quadrature-free technique for the evaluation of integrals with polynomial integrands over B-Reps enclosed by trimmed non-rational spline surfaces.

This technique is based on two successive applications of the divergence theorem, transforming 3D integrals into line integrals that are eventually computed analytically.
The involved steps require the creation and manipulation of (potentially) very high-degree polynomials. Nevertheless, we do not perform explicit evaluation of such functions, but just operations as additions or multiplications (using Bernstein bases), that are known to be more stable.
The accuracy of this integration method has been verified numerically by evaluating integrals of low order polynomials over 2D and 3D domains and comparing the obtained results against reference solutions computed through boundary-conformal quadrature schemes.

In order to apply such integration method to the resolution of PDEs over CAD models using immersed Galerkin discretizations, we transform the integrands of the finite element operators into polynomials.
Thus, relying on~\cite{Mantzaflaris_2015} we create local polynomial approximations of those integrands for every element.
In addition, according to~\cite{Antolin_2019b}, we also approximate at element level the rational B-splines, that may define the geometry, as non-rational B\'ezier curves and surfaces.
This opens the door to the application of the method to B-Reps enclosed by rational splines.

The combination of the results in~\cite{Mantzaflaris_2015,Antolin_2019b} theoretically guarantees the optimal approximation properties of the proposed method for elliptic problems.
This approach is directly extendable to other non-elliptic problems, however suitable approximation properties are not backed up by theoretical evidences.

A series of numerical experiments support our claims.
Thus, the method's performance is illustrated by a series of elliptic problems on immersed 2D and 3D geometries, some of which present rational geometries.
Optimal convergence rates were confirmed in all the cases.
Finally, and in order to prove the potential of the method, its real applicability is demonstrated with a couple of 3D B-Rep models with an industrial level of geometrical complexity.

In this work, we particularize our study to the case of isogeometric discretizations.
Nevertheless, the ideas behind are straightforwardly extendable to other immersed methods as, for instance, the finite cell method or CutFEM/IGA~\cite{wassermann2017geometric,burman2015cutfem,elfverson2018cutiga}, or to other discretization techniques like XFEM or high-order virtual element~\cite{sukumar2000extended,beirao2013basic} methods.
In addition, the quadrature-free integration could be also handy for the evaluation of the right-hand-side integrals involved in moment fitting techniques~\cite{Mueller_2013}.

\section*{Acknowledgments} 
P.\ Antolin was partially supported by the European Research Council through the H2020 ERC Advanced Grant 2015 n.694515 CHANGE,
and by the Swiss National Science Foundation through the project ``Design-through-Analysis (of PDEs): the litmus test'' n.40B2-0\_187094 (BRIDGE Discovery 2019).
T.\ Hirschler was supported by the European Union's Horizon 2020 research and innovation program under grant agreement n.862025 (ADAM2).

\appendix
{\normalsize
  \section{Bernstein polynomials} \label{sec:appA}
In this Appendix we discuss the construction of polynomials using Bernstein bases.
We first introduce, in~\ref{sec:bernstein}, the Bernstein basis, its main properties, and the construction of univariate polynomials.
Afterwards, in~\ref{sec:multi:bernstein}, we discuss its generalization to the case of tensor-product polynomials.
And finally, in~\ref{sec:vector:bernstein} we present the case of multi-dimensional vector polynomials.
Most of the constructions detailed in this Appendix are rather classical and can be found, for instance, in~\cite{Farouki_1988}.

\subsection{Bernstein basis and univariate polynomials} \label{sec:bernstein}
Let us first introduce the Bernstein polynomials basis for a degree $p\geq0$:
\begin{equation}
    B_{i}^p(t) = \binom{p}{i} t^i(1-t)^{p-i}, ~~i=0,\dots,p.
\end{equation}
It is well-known that this basis constitutes an appealing alternative to monomials in terms of numerical stability when it comes to floating-point operations~\cite{Farouki_1988}.

In addition, the Bernstein basis presents some handy properties that simplify the manipulation of polynomials.
For instance, their derivatives can be easily computed as a function of lower degree polynomials. Thus, for $p>0$:
\begin{subequations} \label{eq:berns:der}
  \begin{align}
  \frac{\diff B_{0}^p(t)}{\diff t} &= - p B_{0}^{p-1}(t)\,,\\
  \frac{\diff B_{i}^p(t)}{\diff t} &= p B_{i-1}^{p-1}(t) - p B_{i}^{p-1}(t)\quad\text{for }0<i<p\,,\\
  \frac{\diff B_{p}^p(t)}{\diff t} &= p B_{p-1}^{p-1}(t)\,.
  \end{align}
\end{subequations}
In the same way, their primitives can be computed using polynomials of higher degree:
\begin{equation} \label{eq:berns:anti}
  \int B_{i}^p(t)\diff t = \frac{1}{p+1}\sum_{j=i+1}^{p+1} B_{j}^{p+1}(t)\,,
\end{equation}
that yields:
\begin{equation} \label{eq:berns:intg}
  \int_{0}^{1} B_{i}^p(t)\diff t = \frac{1}{p+1}\,.
\end{equation}
In addition, any Bernstein polynomial of degree $p-1$, with $p>0$, can be expressed as a combination of polynomials of degree $p$ as:
\begin{equation} \label{eq:berns:rise}
B_{i}^{p-1}(t) = \frac{p-i}{p} B_{i}^{p}(t) + \frac{i+1}{p} B_{i+1}^{p}(t)\,.
\end{equation}

Using the Bernstein basis, a univariate real polynomial $f(t)$ of degree $p$ can be expressed as:
\begin{equation} \label{eq:pol:f}
f(t) = \sum_{i=0}^p B_{i}^p(t) f_i,
\end{equation}
where $f_i\in\mathbb{R}$.
Applying~\eqref{eq:berns:der}, \eqref{eq:berns:anti}, \eqref{eq:berns:intg}, and \eqref{eq:berns:rise} to each Bernstein basis function of the polynomial $f(t)$, it is straightforward to compute the derivative of $f(t)$, its antiderivative, integrate it over the domain $[0,1]$, and express it using a basis of degree $p+1$, respectively.
In particular, due to its particular interest in this work, the integral of $f(t)$ over the domain $[0,1]$ is detailed:
\begin{equation} \label{eq:pol:intg}
  \int_{0}^{1} f(t)\diff t = \frac{1}{p+1}\sum_{i=0}^{p}f_i\,.
\end{equation}
This result can be directly applied to the computation of the integral~\eqref{eq:int_1D_4}, in Section~\ref{sec:integration}.
We remark that in this operation no polynomial evaluations are involved, simply the linear combination of the coefficients $f_i$, what makes this computation stable even for high degree polynomials.

Let us know introduce now a second polynomial $g(t)$ of degree $q\geq 0$:
\begin{equation}
g(t) = \sum_{i=0}^q B_{i}^q(t) g_i\,.
\end{equation}
In the case $q=p$, the addition (subtraction) of $f(t)$ and $g(t)$ it is easily computed by adding (subtracting) their coefficients:
\begin{equation} \label{eq:pol:sum}
f(t) \pm g(t) = \sum_{i=0}^p B_{i}^q(t) (f_i \pm g_i)\,.
\end{equation}
On the other hand, if $q < p$, $g(t)$ must be firstly written in the Bernstein basis of degree $p$, applying~\eqref{eq:berns:rise} $p-q$ times, and then the expression~\eqref{eq:pol:sum} can be directly used.

The multiplication of polynomials is another operation that is extensively used in Section~\ref{sec:integration}.
The product $f(t) g(t)$ yields a new polynomial of degree $p+q$ that can be computed as:
\begin{equation} \label{eq:pol:prod}
  f(t)\,g(t) = \sum_{i=0}^{p+q}B_{i}^{p+q}(t)\Bigg(\sum_{j=\max(0,i-q)}^{\min(p,i)}\frac{\binom{p}{j}\binom{q}{i-j}}{\binom{p+q}{i}} f_j g_{i-j}\Bigg).
\end{equation}
Based on that, the composition of two polynomials $f\circ g(t)$ is expressed as:
\begin{equation} \label{eq:berns:composition}
  \begin{split}
    f\circ g(t)
    &= \sum_{i=0}^{p} B_{i}^{p}\left( g(t)\right) f_i\\
    &= \sum_{i=0}^{p} \binom{p}{i} g(t)^i \left(1-g(t)\right)^{p-i}f_i,
  \end{split}
\end{equation}
where the terms $g(t)^i \left(1-g(t)\right)^{p-i}$, $i=1,\dots,p$, can be evaluated by means of the polynomials product expression~\eqref{eq:pol:prod}.

\subsection{Multivariate polynomials} \label{sec:multi:bernstein}
The univariate construction~\eqref{eq:pol:f} can be extended to the case of $m$-dimensional tensor-product polynomials as:
\begin{multline}
  \label{eq:multipol:h}
  h(t_1,t_2,\dots,t_m)  = \\
  \sum_{i_1=0}^{p_1}\sum_{i_2=0}^{p_2}\cdots\sum_{i_m=0}^{p_m}
  B_{i_1}^{p_1}(t_1) B_{i_2}^{p_2}(t_2) \cdots B_{i_m}^{p_m}(t_m) h_{\bm{i}}\,,
\end{multline}
where $(p_1,p_2,\cdots,p_m)$ are the non-negative degrees along the $m$ parametric directions, and $\bm{i}=(i_1,i_2,\cdots,i_m)$ the multi-index accounting for all the univariate indices.
Operations defined for univariate polynomials, as derivatives~\eqref{eq:berns:der}, primitives~\eqref{eq:berns:anti}, or degree raising~\eqref{eq:berns:rise}, can now be applied for every parametric direction independently.
For instance, the computation of antiderivatives along different directions is required in Equations~\eqref{eq:antiderivative_0} and  \eqref{eq:antiderivative_1}.
On the other hand, the integral of $h(t_1,\dots,t_m)$ over a domain $[0,1]^m$ can be easily computed as:
\begin{multline} \label{eq:multipol:intg}
  \int_{0}^{1}\int_{0}^{1}\dots\int_{0}^{1} h(t_1, t_2,\dots,t_m)\diff t_1\diff t_2\dots\diff t_m \\=
  \frac{1}{\prod^{m}_{j=1} \left(p_j+1\right)} \sum_{\bm{i}}h_{\bm{i}}\,.
\end{multline}
As discussed in Remark~\ref{remark:nontrimmed2D}, for the case of non-trimmed B\'ezier patches, the surface integral~\eqref{eq:int_2D_1} can be directly computed using the expression above.
The same applies to the integral~\eqref{eq:goalIntegral} in the case the integration domain is a unit cube (what is applicable to integration over non-cut elements $Q\in\T^{\text{int}}_h(\Omega)$ as discussed in Section~\ref{sec:immersed_iga}).
This is the case of the evaluation of integrals over the non-cut elements discussed in Section~\ref{sec:immersed_iga}.
We also remark here that, as for the univariate case~\eqref{eq:berns:intg}, no polynomial evaluations are required for computing this integral, only a linear combination of the coefficients $h_{\bm{i}}$.

We now consider a second $m$-dimensional polynomial $l(t_1,\dots,t_m)$ with non-negative degrees $(q_1,\dots,q_m)$:
\begin{multline} \label{eq:multipol:l}
l(t_1,t_2,\dots,t_m)  = \\\sum_{i_1=0}^{q_1}\sum_{i_2=0}^{q_2}\cdots\sum_{i_m=0}^{q_m}
B_{i_1}^{q_1}(t_1) B_{i_2}^{q_2}(t_2) \cdots B_{i_m}^{q_m}(t_m) l_{\bm{i}}\,,
\end{multline}
where $l_{\bm{i}}\in\mathbb{R}$.
The multiplication of two $m$-dimensional polynomials, analogously to~\eqref{eq:pol:prod}, results in a polynomial with degrees $(p_1+q_1,\dots,p_m+q_m)$ that can be computed as:
{
  \begin{strip}
    \begin{multline} \label{eq:multipol:prod}
      h(t_,\dots,t_m)\,l(t_,\dots,t_m) =
      \sum_{i_1=0}^{p_1+q_1}
      \cdots
      \sum_{i_m=0}^{p_m+q_m}
      B_{i_1}^{p_1+q_1}(t_1)\cdots B_{i_m}^{p_m+q_m}(t_m)\\
      \sum_{j_1=\max(0,i_1-q_1)}^{\min(p_1,i_1)}
      \cdots
      \sum_{j_m=\max(0,i_m-q_m)}^{\min(p_m,i_m)}
      \frac{\binom{p_1}{j_1}\binom{q_1}{i_1-j_1}}{\binom{p_1+q_1}{i_1}}
      \cdots
      \frac{\binom{p_m}{j_m}\binom{q_m}{i_m-j_m}}{\binom{p_m+q_m}{i_m}}
      h_{j_1,\dots,j_m} l_{i_1-j_1,\dots,i_m-j_m}.
    \end{multline}%
  \end{strip}
}

\subsection{Vector polynomials} \label{sec:vector:bernstein}
The univariate and multivariate polynomials studied above constitute the foundation for the construction of B\'ezier curves, surfaces, and other higher dimensional geometric objects.
In particular, following the polynomial constructions~\eqref{eq:pol:f} and~\eqref{eq:multipol:h}, B\'ezier curves and surfaces can be expressed as:
\begin{subequations} 
  \begin{align}
  \bm{c}(t) &= \sum_{i=0}^p B_{i}^p(t) \bm{c}_i,\\
  \bm{S}(t_1, t_2) &= \sum_{i_1=0}^{p_1}\sum_{i_2=0}^{p_2} B_{i_1}^{p_1}(t_1)B_{i_2}^{p_2}(t_2) \bm{S}_{i_1,i_2},
  \end{align}
\end{subequations}
where $\bm{c}_i,\bm{S}_{i_1,i_2}\in\mathbb{R}^d$ and $d$ is the spatial dimension.
The single coordinate components of $\bm{c}$ and $\bm{S}$ are themselves scalar polynomials and can expressed as:
\begin{subequations} \label{eq:bezier:comp}
  \begin{align}
    c_k(t)  &= \bm{c}(t)\cdot\bm{e}_k = \sum_{i=0}^p B_{i}^p(t) \bm{c}_i\cdot\bm{e}_k,\\
    S_k(t_1,t_2)
            &= \bm{S}(t_1,t_2)\cdot\bm{e}_k\nonumber \\
            &= \sum_{i_1=0}^{p_1} \sum_{i_2=0}^{p_2} B_{i_1}^{p_1}(t) B_{i_2}^{p_2}(t) \bm{S}_{i_1,i_2}\cdot\bm{e}_k,
  \end{align}
\end{subequations}
for $k=1,\dots,d$, and where $\bm{e}_k$ are the unit vectors along the Cartesian directions.

Thus, operations like partial derivatives, or cross and scalar products between B\'ezier curves and surfaces, like the ones used in Section~\ref{sec:integration}, can be carried out by using its individual coordinate components~\eqref{eq:bezier:comp} and combining them according to the operations detailed in previous sections for scalar univariate and multivariate polynomials.
Among all the operations, due to its higher complexity, in what remains we detail the composition between multivariate B\'eziers.

We consider two multivariate B\'eziers $\bm{F}:\mathbb{R}^s\to\mathbb{R}^d$ and $\bm{G}:\mathbb{R}^m\to\mathbb{R}^s$ of the form:
{\small
  \begin{subequations}
    \begin{align}
      \bm{F}(r_1,\dots,r_m)  &= \sum_{i_1=0}^{p_1}\cdots\sum_{i_s=0}^{p_s} B_{i_1}^{p_1}(t_1) \cdots B_{i_s}^{p_s}(t_s) \bm{F}_{\bm{i}}\,,\\
      \bm{G}(t_1,\dots,t_s)  &= \sum_{j_1=0}^{q_1}\cdots\sum_{j_m=0}^{q_m} B_{j_1}^{q_1}(r_1) \cdots B_{j_m}^{q_m}(r_m) \bm{G}_{\bm{j}}\,,
    \end{align}
  \end{subequations}%
}
that have non-negative degrees $(p_1,p_2,\dots,p_s)$ and $(q_1,q_2,\dots,q_m)$, respectively.
$\bm{G}_{\bm{i}}\in\mathbb{R}^s$ and $\bm{F}_{\bm{i}}\in\mathbb{R}^d$ are the associated control points, and $\mathbf{i}=(i_1,i_2,\dots,i_s)$ and $\mathbf{j}=(j_1,j_2,\dots,j_m)$ the corresponding multi-indices.
We want to compute the composition $\bm{F}\circ\bm{G} (t_1,\dots,t_m):\mathbb{R}^m\to\mathbb{R}^d$.
Working with the coordinate components $G_k(t_1,\dots,t_m) = \bm{G}(t_1,\dots,t_m)\cdot\bm{e}_k,~k=1,\dots,s$, we obtain:
{\small
  \begin{multline}
      \bm{F}\circ\bm{G}(t_1,\dots,t_m)  = \sum_{i_1=0}^{p_1}\cdots\sum_{i_s=0}^{p_s}
      B_{i_1}^{p_1}\circ G_1(t_1,\dots,t_m) \\\cdots B_{i_s}^{p_s}\circ G_s(t_1,\dots,t_m) \bm{F}_{\bm{i}}\,.
  \end{multline}%
}
Every term $B_{i_k}^{q_k}\circ G_k(t_1,\dots,t_m),~k=1,\dots,s$, is the composition between a univariate Bernstein polynomial and a $m$-dimensional scalar polynomial expressed in a tensor-product Bernstein basis:%
{\small
  \begin{multline}
    B_{i_k}^{p_k}\circ G_k(t_1,\dots,t_m) = \\
    \binom{p_k}{i_k} G_k\left(t_1,\dots,t_m\right)^{i_k}\left(1-G_k\left(t_1,\dots,t_m\right)\right)^{p_k-i_k},
  \end{multline}%
}
where the products
are computed performing multiplications between multi-dimensional scalar polynomials, detailed in Equation~\eqref{eq:multipol:prod}.
%
}

\bibliography{biblio}

\begin{thebibliography}{10}
\expandafter\ifx\csname url\endcsname\relax
  \def\url#1{\texttt{#1}}\fi
\expandafter\ifx\csname urlprefix\endcsname\relax\def\urlprefix{URL }\fi
\expandafter\ifx\csname href\endcsname\relax
  \def\href#1#2{#2} \def\path#1{#1}\fi

\bibitem{Hughes_2005}
T.~J.~R. Hughes, J.~A. Cottrell, Y.~Bazilevs, Isogeometric analysis: {CAD},
  finite elements, {NURBS}, exact geometry and mesh refinement, Computer
  Methods in Applied Mechanics and Engineering 194~(39-41) (2005) 4135--4195.
\newblock \href {https://doi.org/10.1016/j.cma.2004.10.008}
  {\path{doi:10.1016/j.cma.2004.10.008}}.

\bibitem{Cottrell_2009}
J.~A. Cottrell, T.~J.~R. Hughes, Y.~Bazilevs, Isogeometric Analysis, John Wiley
  \& Sons, 2009.

\bibitem{Liu_2009}
G.~Liu, Meshfree Methods, {CRC} Press, 2009.
\newblock \href {https://doi.org/10.1201/9781420082104}
  {\path{doi:10.1201/9781420082104}}.

\bibitem{BAZILEVS_2006}
Y.~Bazilevs, L.~B. da~Veiga, J.~A. Cottrell, T.~J.~R. Hughes, G.~Sangalli,
  {Isogeometric} {Analysis}: {approximation}, {stability} {and} {error}
  {estimates} {for} h-{refined} {meshes}, Mathematical Models and Methods in
  Applied Sciences 16~(07) (2006) 1031--1090.
\newblock \href {https://doi.org/10.1142/s0218202506001455}
  {\path{doi:10.1142/s0218202506001455}}.

\bibitem{Buffa_2011}
A.~Buffa, J.~Rivas, G.~Sangalli, R.~V{\'{a}}zquez, Isogeometric discrete
  differential forms in three dimensions, {SIAM} Journal on Numerical Analysis
  49~(2) (2011) 818--844.
\newblock \href {https://doi.org/10.1137/100786708}
  {\path{doi:10.1137/100786708}}.

\bibitem{Hiemstra_2014}
R.~Hiemstra, D.~Toshniwal, R.~Huijsmans, M.~Gerritsma, High order geometric
  methods with exact conservation properties, Journal of Computational Physics
  257 (2014) 1444--1471.
\newblock \href {https://doi.org/10.1016/j.jcp.2013.09.027}
  {\path{doi:10.1016/j.jcp.2013.09.027}}.

\bibitem{Lipton_2010}
S.~Lipton, J.~Evans, Y.~Bazilevs, T.~Elguedj, T.~J.~R. Hughes, Robustness of
  isogeometric structural discretizations under severe mesh distortion,
  Computer Methods in Applied Mechanics and Engineering 199~(5-8) (2010)
  357--373.
\newblock \href {https://doi.org/10.1016/j.cma.2009.01.022}
  {\path{doi:10.1016/j.cma.2009.01.022}}.

\bibitem{Herrema_2017}
A.~J. Herrema, N.~M. Wiese, C.~N. Darling, B.~Ganapathysubramanian,
  A.~Krishnamurthy, M.-C. Hsu, A framework for parametric design optimization
  using isogeometric~analysis, Computer Methods in Applied Mechanics and
  Engineering 316 (2017) 944--965.
\newblock \href {https://doi.org/10.1016/j.cma.2016.10.048}
  {\path{doi:10.1016/j.cma.2016.10.048}}.

\bibitem{Antolin_2019}
P.~Antolin, A.~Buffa, E.~Cohen, J.~F. Dannenhoffer, G.~Elber, S.~Elgeti,
  R.~Haimes, R.~Riesenfeld, Optimizing micro-tiles in micro-structures as a
  design paradigm, Computer-Aided Design 115 (2019) 23--33.
\newblock \href {https://doi.org/10.1016/j.cad.2019.05.020}
  {\path{doi:10.1016/j.cad.2019.05.020}}.

\bibitem{Hafner_2019}
C.~Hafner, C.~Schumacher, E.~Knoop, T.~Auzinger, B.~Bickel, M.~Bächer, {X-CAD:
  Optimizing {CAD} Models with Extended Finite Elements}, {ACM} Transactions on
  Graphics 38~(6) (2019) 1--15.
\newblock \href {https://doi.org/10.1145/3355089.3356576}
  {\path{doi:10.1145/3355089.3356576}}.

\bibitem{Hirschler_2020}
T.~Hirschler, R.~Bouclier, A.~Duval, T.~Elguedj, J.~Morlier, A new lighting on
  analytical discrete sensitivities in the context of {IsoGeometric} shape
  optimization, Archives of Computational Methods in Engineering 28~(4) (2020)
  2371--2408.
\newblock \href {https://doi.org/10.1007/s11831-020-09458-6}
  {\path{doi:10.1007/s11831-020-09458-6}}.

\bibitem{wang2012converting}
W.~Wang, Y.~Zhang, G.~Xu, T.~J.~R. Hughes, Converting an unstructured
  quadrilateral/hexahedral mesh to a rational {T}-spline, Computational
  Mechanics 50~(1) (2012) 65--84.

\bibitem{wei2018blended}
X.~Wei, Y.~J. Zhang, D.~Toshniwal, H.~Speleers, X.~Li, C.~Manni, J.~A. Evans,
  T.~J.~R. Hughes, Blended {B}-spline construction on unstructured
  quadrilateral and hexahedral meshes with optimal convergence rates in
  isogeometric analysis, Computer Methods in Applied Mechanics and Engineering
  341 (2018) 609--639.

\bibitem{xia2017isogeometric}
S.~Xia, X.~Qian, Isogeometric analysis with {B}{\'e}zier tetrahedra, Computer
  Methods in Applied Mechanics and Engineering 316 (2017) 782--816.

\bibitem{Al_Akhras_2016}
H.~A. Akhras, T.~Elguedj, A.~Gravouil, M.~Rochette, Isogeometric
  analysis-suitable trivariate {NURBS} models from standard {B-Rep} models,
  Computer Methods in Applied Mechanics and Engineering 307 (2016) 256--274.
\newblock \href {https://doi.org/10.1016/j.cma.2016.04.028}
  {\path{doi:10.1016/j.cma.2016.04.028}}.

\bibitem{Hinz_2018}
J.~Hinz, M.~M{\"o}ller, C.~Vuik, Elliptic grid generation techniques in the
  framework of isogeometric analysis applications, Computer Aided Geometric
  Design 65 (2018) 48--75.
\newblock \href {https://doi.org/10.1016/j.cagd.2018.03.023}
  {\path{doi:10.1016/j.cagd.2018.03.023}}.

\bibitem{Massarwi_2019}
F.~Massarwi, P.~Antolin, G.~Elber, Volumetric untrimming: Precise decomposition
  of trimmed trivariates into tensor products, Computer Aided Geometric Design
  71 (2019) 1--15.
\newblock \href {https://doi.org/10.1016/j.cagd.2019.04.005}
  {\path{doi:10.1016/j.cagd.2019.04.005}}.

\bibitem{Maquart_2020}
T.~Maquart, Y.~Wenfeng, T.~Elguedj, A.~Gravouil, M.~Rochette, {3D} volumetric
  isotopological meshing for finite element and isogeometric based reduced
  order modeling, Computer Methods in Applied Mechanics and Engineering 362
  (2020) 112809.
\newblock \href {https://doi.org/10.1016/j.cma.2019.112809}
  {\path{doi:10.1016/j.cma.2019.112809}}.

\bibitem{rank2012geometric}
E.~Rank, M.~Ruess, S.~Kollmannsberger, D.~Schillinger, A.~D{\"u}ster, Geometric
  modeling, isogeometric analysis and the finite cell method, Computer Methods
  in Applied Mechanics and Engineering 249 (2012) 104--115.

\bibitem{Legrain_2013}
G.~Legrain, A {NURBS} enhanced extended finite element approach for unfitted
  {CAD} analysis, Computational Mechanics 52~(4) (2013) 913--929.
\newblock \href {https://doi.org/10.1007/s00466-013-0854-7}
  {\path{doi:10.1007/s00466-013-0854-7}}.

\bibitem{Breitenberger_2015}
M.~Breitenberger, A.~Apostolatos, B.~Philipp, R.~Wüchner, K.-U. Bletzinger,
  Analysis in computer aided design: Nonlinear isogeometric {B-Rep} analysis of
  shell structures, Computer Methods in Applied Mechanics and Engineering 284
  (2015) 401--457.
\newblock \href {https://doi.org/10.1016/j.cma.2014.09.033}
  {\path{doi:10.1016/j.cma.2014.09.033}}.

\bibitem{Hsu_2016}
M.-C. Hsu, C.~Wang, F.~Xu, A.~J. Herrema, A.~Krishnamurthy, Direct
  immersogeometric fluid flow analysis using {B-rep} {CAD} models, Computer
  Aided Geometric Design 43 (2016) 143--158.
\newblock \href {https://doi.org/10.1016/j.cagd.2016.02.007}
  {\path{doi:10.1016/j.cagd.2016.02.007}}.

\bibitem{guo2018variationally}
Y.~Guo, J.~Heller, T.~J.~R. Hughes, M.~Ruess, D.~Schillinger, Variationally
  consistent isogeometric analysis of trimmed thin shells at finite
  deformations, based on the step exchange format, Computer Methods in Applied
  Mechanics and Engineering 336 (2018) 39--79.

\bibitem{wassermann2019integrating}
B.~Wassermann, S.~Kollmannsberger, S.~Yin, L.~Kudela, E.~Rank, Integrating
  {CAD} and numerical analysis: 'dirty geometry' handling using the finite cell
  method, Computer Methods in Applied Mechanics and Engineering 351 (2019)
  808--835.

\bibitem{Marussig_2017}
B.~Marussig, T.~J.~R. Hughes, A review of trimming in isogeometric analysis:
  Challenges, data exchange and simulation aspects, Archives of Computational
  Methods in Engineering 25~(4) (2017) 1059--1127.
\newblock \href {https://doi.org/10.1007/s11831-017-9220-9}
  {\path{doi:10.1007/s11831-017-9220-9}}.

\bibitem{peskin2002immersed}
C.~S. Peskin, The immersed boundary method, Acta numerica 11 (2002) 479--517.

\bibitem{D_ster_2008}
A.~D{\"u}ster, J.~Parvizian, Z.~Yang, E.~Rank, The finite cell method for
  three-dimensional problems of solid mechanics, Computer Methods in Applied
  Mechanics and Engineering 197~(45-48) (2008) 3768--3782.
\newblock \href {https://doi.org/10.1016/j.cma.2008.02.036}
  {\path{doi:10.1016/j.cma.2008.02.036}}.

\bibitem{schillinger2012isogeometric}
D.~Schillinger, L.~Dede, M.~A. Scott, J.~A. Evans, M.~J. Borden, E.~Rank,
  T.~J.~R. Hughes, An isogeometric design-through-analysis methodology based on
  adaptive hierarchical refinement of {NURBS}, immersed boundary methods, and
  {T}-spline {CAD} surfaces, Computer Methods in Applied Mechanics and
  Engineering 249 (2012) 116--150.

\bibitem{burman2015cutfem}
E.~Burman, S.~Claus, P.~Hansbo, M.~G. Larson, A.~Massing, {CutFEM}:
  discretizing geometry and partial differential equations, International
  Journal for Numerical Methods in Engineering 104~(7) (2015) 472--501.

\bibitem{wassermann2017geometric}
B.~Wassermann, S.~Kollmannsberger, T.~Bog, E.~Rank, From geometric design to
  numerical analysis: a direct approach using the finite cell method on
  constructive solid geometry, Computers \& Mathematics with Applications
  74~(7) (2017) 1703--1726.

\bibitem{elfverson2018cutiga}
D.~Elfverson, M.~G. Larson, K.~Larsson, {CutIGA} with basis function removal,
  Advanced Modeling and Simulation in Engineering Sciences 5~(1) (2018) 1--19.

\bibitem{shephard1991automatic}
M.~S. Shephard, M.~K. Georges, Automatic three-dimensional mesh generation by
  the finite octree technique, International Journal for Numerical methods in
  engineering 32~(4) (1991) 709--749.

\bibitem{abedian2013performance}
A.~Abedian, J.~Parvizian, A.~D{\"u}ster, H.~Khademyzadeh, E.~Rank, Performance
  of different integration schemes in facing discontinuities in the finite cell
  method, International Journal of Computational Methods 10~(03) (2013)
  1350002.

\bibitem{kudela2016smart}
L.~Kudela, N.~Zander, S.~Kollmannsberger, E.~Rank, Smart octrees: Accurately
  integrating discontinuous functions in {3D}, Computer Methods in Applied
  Mechanics and Engineering 306 (2016) 406--426.

\bibitem{peto2020enhanced}
M.~Pet{\"o}, F.~Duvigneau, S.~Eisentr{\"a}ger, Enhanced numerical integration
  scheme based on image-compression techniques: application to fictitious
  domain methods, Advanced Modeling and Simulation in Engineering Sciences 7
  (2020) 1--42.

\bibitem{Verhoosel_2015}
C.~Verhoosel, G.~van Zwieten, B.~van Rietbergen, R.~de~Borst, Image-based
  goal-oriented adaptive isogeometric analysis with application to the
  micro-mechanical modeling of trabecular bone, Computer Methods in Applied
  Mechanics and Engineering 284 (2015) 138--164.
\newblock \href {https://doi.org/10.1016/j.cma.2014.07.009}
  {\path{doi:10.1016/j.cma.2014.07.009}}.

\bibitem{Divi_2020}
S.~C. Divi, C.~V. Verhoosel, F.~Auricchio, A.~Reali, E.~H. van Brummelen,
  Error-estimate-based adaptive integration for immersed isogeometric analysis,
  Computers {\&} Mathematics with Applications 80~(11) (2020) 2481--2516.
\newblock \href {https://doi.org/10.1016/j.camwa.2020.03.026}
  {\path{doi:10.1016/j.camwa.2020.03.026}}.

\bibitem{Kudela_2015}
L.~Kudela, N.~Zander, T.~Bog, S.~Kollmannsberger, E.~Rank, Efficient and
  accurate numerical quadrature for immersed boundary methods, Advanced
  Modeling and Simulation in Engineering Sciences 2~(1) (2015).
\newblock \href {https://doi.org/10.1186/s40323-015-0031-y}
  {\path{doi:10.1186/s40323-015-0031-y}}.

\bibitem{Antolin_2019b}
P.~Antolin, A.~Buffa, M.~Martinelli, {Isogeometric Analysis on {V-reps}: First
  results}, Computer Methods in Applied Mechanics and Engineering 355 (2019)
  976--1002.
\newblock \href {https://doi.org/10.1016/j.cma.2019.07.015}
  {\path{doi:10.1016/j.cma.2019.07.015}}.

\bibitem{joulaian2016numerical}
M.~Joulaian, S.~Hubrich, A.~D{\"u}ster, Numerical integration of
  discontinuities on arbitrary domains based on moment fitting, Computational
  Mechanics 57~(6) (2016) 979--999.

\bibitem{Hubrich_2017}
S.~Hubrich, P.~D. Stolfo, L.~Kudela, S.~Kollmannsberger, E.~Rank, A.~Schröder,
  A.~D{\"u}ster, Numerical integration of discontinuous functions: moment
  fitting and smart octree, Computational Mechanics 60~(5) (2017) 863--881.
\newblock \href {https://doi.org/10.1007/s00466-017-1441-0}
  {\path{doi:10.1007/s00466-017-1441-0}}.

\bibitem{Hubrich_2019}
S.~Hubrich, A.~D{\"u}ster, Numerical integration for nonlinear problems of the
  finite cell method using an adaptive scheme based on moment fitting,
  Computers {\&} Mathematics with Applications 77~(7) (2019) 1983--1997.
\newblock \href {https://doi.org/10.1016/j.camwa.2018.11.030}
  {\path{doi:10.1016/j.camwa.2018.11.030}}.

\bibitem{bui2020efficient}
H.-G. Bui, D.~Schillinger, G.~Meschke, Efficient cut-cell quadrature based on
  moment fitting for materially nonlinear analysis, Computer Methods in Applied
  Mechanics and Engineering 366 (2020) 113050.

\bibitem{Lasserre_1998}
J.~B. Lasserre, Integration on a convex polytope, Proceedings of the American
  Mathematical Society 126~(8) (1998) 2433--2441.
\newblock \href {https://doi.org/10.1090/s0002-9939-98-04454-2}
  {\path{doi:10.1090/s0002-9939-98-04454-2}}.

\bibitem{Gonzalez_Ochoa_1998}
C.~Gonzalez-Ochoa, S.~McCammon, J.~Peters, Computing moments of objects
  enclosed by piecewise polynomial surfaces, {ACM} Transactions on Graphics
  17~(3) (1998) 143--157.
\newblock \href {https://doi.org/10.1145/285857.285858}
  {\path{doi:10.1145/285857.285858}}.

\bibitem{Mousavi_2010}
S.~E. Mousavi, N.~Sukumar, Numerical integration of polynomials and
  discontinuous functions on irregular convex polygons and polyhedrons,
  Computational Mechanics 47~(5) (2010) 535--554.
\newblock \href {https://doi.org/10.1007/s00466-010-0562-5}
  {\path{doi:10.1007/s00466-010-0562-5}}.

\bibitem{Chin_2015}
E.~B. Chin, J.~B. Lasserre, N.~Sukumar, Numerical integration of homogeneous
  functions on convex and nonconvex polygons and polyhedra, Computational
  Mechanics 56~(6) (2015) 967--981.
\newblock \href {https://doi.org/10.1007/s00466-015-1213-7}
  {\path{doi:10.1007/s00466-015-1213-7}}.

\bibitem{Chin_2020}
E.~B. Chin, N.~Sukumar, An efficient method to integrate polynomials over
  polytopes and curved solids, Computer Aided Geometric Design 82 (2020)
  101914.
\newblock \href {https://doi.org/10.1016/j.cagd.2020.101914}
  {\path{doi:10.1016/j.cagd.2020.101914}}.

\bibitem{Ventura_2006}
G.~Ventura, On the elimination of quadrature subcells for discontinuous
  functions in the {eXtended} {Finite-Element Method}, International Journal
  for Numerical Methods in Engineering 66~(5) (2006) 761--795.
\newblock \href {https://doi.org/10.1002/nme.1570}
  {\path{doi:10.1002/nme.1570}}.

\bibitem{Duczek_2015}
S.~Duczek, U.~Gabbert, Efficient integration method for fictitious domain
  approaches, Computational Mechanics 56~(4) (2015) 725--738.
\newblock \href {https://doi.org/10.1007/s00466-015-1197-3}
  {\path{doi:10.1007/s00466-015-1197-3}}.

\bibitem{Abedian_2019}
A.~Abedian, A.~D{\"u}ster, Equivalent {L}egendre polynomials: Numerical
  integration of discontinuous functions in the finite element methods,
  Computer Methods in Applied Mechanics and Engineering 343 (2019) 690--720.
\newblock \href {https://doi.org/10.1016/j.cma.2018.08.002}
  {\path{doi:10.1016/j.cma.2018.08.002}}.

\bibitem{Mueller_2013}
B.~M{\"u}ller, F.~Kummer, M.~Oberlack, Highly accurate surface and volume
  integration on implicit domains by means of moment-fitting, International
  Journal for Numerical Methods in Engineering 96~(8) (2013) 512--528.
\newblock \href {https://doi.org/10.1002/nme.4569}
  {\path{doi:10.1002/nme.4569}}.

\bibitem{Sudhakar_2014}
Y.~Sudhakar, J.~M. de~Almeida, W.~A. Wall, An accurate, robust, and
  easy-to-implement method for integration over arbitrary polyhedra:
  Application to embedded interface methods, Journal of Computational Physics
  273 (2014) 393--415.
\newblock \href {https://doi.org/10.1016/j.jcp.2014.05.019}
  {\path{doi:10.1016/j.jcp.2014.05.019}}.

\bibitem{Gunderman_2021}
D.~Gunderman, K.~Weiss, J.~A. Evans, High-accuracy mesh-free quadrature for
  trimmed parametric surfaces and volumes, submitted (Jan. 2021).
\newblock \href {http://arxiv.org/abs/2101.06497} {\path{arXiv:2101.06497}}.

\bibitem{parvizian2007finite}
J.~Parvizian, A.~D{\"u}ster, E.~Rank, Finite cell method, Computational
  Mechanics 41~(1) (2007) 121--133.

\bibitem{giannelli2012thb}
C.~Giannelli, B.~J{\"u}ttler, H.~Speleers, {THB}-splines: {The} truncated basis
  for hierarchical splines, Computer Aided Geometric Design 29~(7) (2012)
  485--498.

\bibitem{bazilevs2010isogeometric}
Y.~Bazilevs, V.~M. Calo, J.~A. Cottrell, J.~A. Evans, T.~J. R.~R. Hughes,
  S.~Lipton, M.~A. Scott, T.~W. Sederberg, Isogeometric analysis using
  {T}-splines, Computer Methods in Applied Mechanics and Engineering 199~(5-8)
  (2010) 229--263.

\bibitem{B_chet_2005}
E.~B{\'{e}}chet, H.~Minnebo, N.~Moës, B.~Burgardt, Improved implementation and
  robustness study of the x-{FEM} for stress analysis around cracks,
  International Journal for Numerical Methods in Engineering 64~(8) (2005)
  1033--1056.
\newblock \href {https://doi.org/10.1002/nme.1386}
  {\path{doi:10.1002/nme.1386}}.

\bibitem{de_Prenter_2017}
F.~de~Prenter, C.~Verhoosel, G.~van Zwieten, E.~van Brummelen, Condition number
  analysis and preconditioning of the finite cell method, Computer Methods in
  Applied Mechanics and Engineering 316 (2017) 297--327.
\newblock \href {https://doi.org/10.1016/j.cma.2016.07.006}
  {\path{doi:10.1016/j.cma.2016.07.006}}.

\bibitem{Buffa_2020}
A.~Buffa, R.~Puppi, R.~V{\'{a}}zquez, A minimal stabilization procedure for
  isogeometric methods on trimmed geometries, {SIAM} Journal on Numerical
  Analysis 58~(5) (2020) 2711--2735.
\newblock \href {https://doi.org/10.1137/19m1244718}
  {\path{doi:10.1137/19m1244718}}.

\bibitem{Hansbo_2002}
A.~Hansbo, P.~Hansbo, An unfitted finite element method, based on {Nitsche}'s
  method, for elliptic interface problems, Computer Methods in Applied
  Mechanics and Engineering 191~(47-48) (2002) 5537--5552.
\newblock \href {https://doi.org/10.1016/s0045-7825(02)00524-8}
  {\path{doi:10.1016/s0045-7825(02)00524-8}}.

\bibitem{Ruess_2013}
M.~Ruess, D.~Schillinger, Y.~Bazilevs, V.~Varduhn, E.~Rank, Weakly enforced
  essential boundary conditions for {NURBS}-embedded and trimmed {NURBS}
  geometries on the basis of the finite cell method, International Journal for
  Numerical Methods in Engineering 95~(10) (2013) 811--846.
\newblock \href {https://doi.org/10.1002/nme.4522}
  {\path{doi:10.1002/nme.4522}}.

\bibitem{Pande_2021}
S.~Pande, P.~Papadopoulos, I.~Babu{\v{s}}ka, A cut-cell finite element method
  for {Poisson}'s equation on arbitrary planar domains, Computer Methods in
  Applied Mechanics and Engineering 383 (2021) 113875.
\newblock \href {https://doi.org/10.1016/j.cma.2021.113875}
  {\path{doi:10.1016/j.cma.2021.113875}}.

\bibitem{Mantzaflaris_2015}
A.~Mantzaflaris, B.~J{\"u}ttler, Integration by interpolation and look-up for
  {Galerkin}-based isogeometric analysis, Computer Methods in Applied Mechanics
  and Engineering 284 (2015) 373--400.
\newblock \href {https://doi.org/10.1016/j.cma.2014.09.014}
  {\path{doi:10.1016/j.cma.2014.09.014}}.

\bibitem{borden2011isogeometric}
M.~J. Borden, M.~A. Scott, J.~A. Evans, T.~J.~R. Hughes, Isogeometric finite
  element data structures based on {B}{\'e}zier extraction of {NURBS},
  International Journal for Numerical Methods in Engineering 87~(1-5) (2011)
  15--47.

\bibitem{d2018multi}
D.~D’Angella, S.~Kollmannsberger, E.~Rank, A.~Reali, Multi-level {B}{\'e}zier
  extraction for hierarchical local refinement of isogeometric analysis,
  Computer Methods in Applied Mechanics and Engineering 328 (2018) 147--174.

\bibitem{scott2011isogeometric}
M.~A. Scott, M.~J. Borden, C.~V. Verhoosel, T.~W. Sederberg, T.~J.~R. Hughes,
  Isogeometric finite element data structures based on {B}{\'e}zier extraction
  of {T}-splines, International Journal for Numerical Methods in Engineering
  88~(2) (2011) 126--156.

\bibitem{Cohen2001}
E.~Cohen, R.~F. Riesenfeld, G.~Elber, Geometric Modeling with Splines, Taylor
  \& Francis Ltd., 2001.

\bibitem{Farin2001}
G.~Farin, Curves and Surfaces for {CAGD}: A Practical Guide, Morgan Kaufmann
  Publ inc, 2001.

\bibitem{Piegl_1997}
L.~Piegl, W.~Tiller, The {NURBS} Book, Springer Berlin Heidelberg, 1997.
\newblock \href {https://doi.org/10.1007/978-3-642-59223-2}
  {\path{doi:10.1007/978-3-642-59223-2}}.

\bibitem{requicha1992solid}
A.~A. Requicha, J.~R. Rossignac, Solid modeling and beyond, IEEE computer
  graphics and applications 12~(5) (1992) 31--44.

\bibitem{braid1973designing}
I.~C. Braid, Designing with volumes, Ph.D. thesis, University of Cambridge
  (1973).

\bibitem{Antonietti_2018}
P.~F. Antonietti, P.~Houston, G.~Pennesi, Fast numerical integration on
  polytopic meshes with applications to discontinuous {Galerkin} finite element
  methods, Journal of Scientific Computing 77~(3) (2018) 1339--1370.
\newblock \href {https://doi.org/10.1007/s10915-018-0802-y}
  {\path{doi:10.1007/s10915-018-0802-y}}.

\bibitem{OpenCASCADE}
O.~C. SAS, {Open CASCADE} 7.3.0, \url{http:///www.opencascade.com} (May, 2018).

\bibitem{Irit}
G.~Elber, {Irit} 11 user's manual, \url{http://www.cs.technion.ac.il/~irit/}
  (2019).

\bibitem{antolin2021overlapping}
P.~Antolin, A.~Buffa, R.~Puppi, X.~Wei, Overlapping multipatch isogeometric
  method with minimal stabilization, SIAM Journal on Scientific Computing
  43~(1) (2021) A330--A354.

\bibitem{sukumar2000extended}
N.~Sukumar, N.~Mo{\"e}s, B.~Moran, T.~Belytschko, Extended finite element
  method for three-dimensional crack modelling, International journal for
  numerical methods in engineering 48~(11) (2000) 1549--1570.

\bibitem{beirao2013basic}
L.~Beir{\~a}o~da Veiga, F.~Brezzi, A.~Cangiani, G.~Manzini, L.~D. Marini,
  A.~Russo, Basic principles of virtual element methods, Mathematical Models
  and Methods in Applied Sciences 23~(01) (2013) 199--214.

\bibitem{Farouki_1988}
R.~Farouki, V.~Rajan, Algorithms for polynomials in bernstein form, Computer
  Aided Geometric Design 5~(1) (1988) 1--26.
\newblock \href {https://doi.org/10.1016/0167-8396(88)90016-7}
  {\path{doi:10.1016/0167-8396(88)90016-7}}.

\end{thebibliography}

\end{document}